\documentclass[amstex,12pt]{article}%
\usepackage{amsfonts}
\usepackage{amsmath}
\usepackage{amssymb}
\usepackage{graphicx}%
\providecommand{\U}[1]{\protect\rule{.1in}{.1in}}

\begin{document}

\title{{\Large Adjoint Difference Equation for the Nikiforov-Uvarov-Suslov Difference
Equation of Hypergeometric Type on Non-uniform Lattices }}
\author{Jinfa Cheng$^{1\ast}$, Weizhong Dai$^{2}$}
\date{{\small 1. School of Mathematical Sciences, Xiamen University, }\\
{\small Xiamen, Fujian, 361005, P. R. China}\\
{\small 2. Mathematics \& Statistics, College of Engineering \& Science,}\\
{\small Louisiana Tech University, Ruston, LA 71272, USA}\\
{\small E-mail: dai@coes.latech.edu}\\
{\small *Corresponding\ author. E-mail: jfcheng@xmu.edu.cn}}
\maketitle

\begin{abstract}
In this article, we obtain the adjoint difference equation for the
Nikiforov-Uvarov-Suslov difference equation of hypergeometric type on
non-uniform lattices, and prove it to be a difference equation of
hypergeometric type on non-uniform lattices as well. The particular solutions
of the adjoint difference equation are then obtained. As an application of
these particular solutions, we use them to obtain the particular solutions for
the original difference equation of hypergeometric type on non-uniform
lattices. In addition, we give another kind of fundamental theorems for the
Nikiforov-Uvarov-Suslov difference equation of hypergeometric type, which are
essentially new results and their expressions are different from the Suslov
Theorem. Finally, we give an example to illustrate the application of the new
fundamental theorems.

\end{abstract}



\section{{\protect\Large Introduction}}

Differential equation of hypergeometric type:%
\begin{equation}
\sigma(z)y^{\prime\prime}(z)+\tau(z)y^{\prime}(z)+\lambda y(z)=0,
\label{hytype}%
\end{equation}
where $\sigma(z)$ and $\tau(z)$ are polynomials of degrees at most two and
one, respectively, and $\lambda$ is a constant, has attracted great attention,
since its solutions are some types of special functions of mathematical
physics, such as the classical orthogonal polynomials, the hypergeometric and
cylindrical functions. In particular, for some positive integer $n$ such that
$\lambda=-\frac{n(n-1)\sigma^{\prime\prime}}{2}-n\tau^{\prime}~$%
and$~\lambda_{m}\neq\lambda_{n}~$for$~m=0,1,\dots,n-1$, Eq. (\ref{hytype}) has
a polynomial solution $y_{n}(z)$ of degree $n$, which can be expressed by the
Rodrigues formula \cite{nikiforov1991, nikiforov1988, suslov1989, wang1989,
andrews1985,andrews1999} as
\begin{equation}
y_{n}(z)=\frac{1}{\rho(z)}\frac{d^{n}}{dz^{n}}(\rho(z)\sigma^{n}(z)),
\end{equation}
where $\rho(z)$ satisfies the Pearson equation
\begin{equation}
(\sigma(z)\rho(z))^{\prime}=\tau(z)\rho(z).
\end{equation}
These solution functions are useful in quantum mechanics, the theory of group
representations, and computational mathematics. Because of this, the classical
theory of hypergeometric type equations has been greatly developed by G.
Andrews, R. Askey \cite{andrews1985, andrews1999}, J.A. Wilson, M. Ismail
\cite{askey1979, askey1984, askey1985,ismail1992a}; F. Nikiforov, K. Suslov,
B. Uvarov, N.M. Atakishiyev \cite{nikiforov1991, nikiforov1988, suslov1989,
nikiforov1983, suslov1992, suslov1988}; G. George, M. Rahman \cite{gasper2004}%
; T.H. Koornwinder \cite{koornwinder1994}; and many other researchers
\cite{Foup2008, Foup2013, jia2017, Magnus1995, swarttouw1994, Witte2011,
dreyfus2015, horner1963, horner1964, kac2002, robin2013,Bangerezako2005}. On
the other hand, many researchers like R. \'{A}lvarez-Nodarse, K.L. Cardoso, I.
Area, E. Godoy, A. Ronveaux, A. Zarzo \cite{alvarez2011, area2003, area2005}
studied particular solutions for the adjoint differential equation of Eq.
(\ref{hytype}) as
\begin{equation}
\sigma(z)w^{\prime\prime}(z)+(2\sigma^{\prime}(z)-\tau(z))w^{\prime
}(z)+(\lambda+\sigma^{\prime\prime}(z)-\tau^{\prime}(z))w(z)=0,
\end{equation}
or the alternative one as
\begin{equation}
\sigma(z)y^{\prime\prime}(z)-\tau_{-2}(z)y^{\prime}(z)+(\lambda-\kappa
_{-1})y(z)=0,
\end{equation}
where
\begin{equation}
\tau_{\nu}(z)\doteq\tau(z)+\nu\sigma^{\prime}(z),\kappa_{\nu}\doteq
\tau^{\prime}+\frac{1}{2}(\nu-1)\sigma^{\prime\prime}. \label{tauka}%
\end{equation}

A.F. Nikiforov, V.B. Uvarov and S.K. Suslov \cite{nikiforov1991,
nikiforov1988} generalized Eq. (\ref{hytype}) to a difference equation of
hypergeometric type case and studied the Nikiforov-Uvarov-Suslov difference
equation on a lattice $x(s)$ with variable step size $\nabla x(s)=x(s)-x(s-1)$
as%
\begin{equation}
\widetilde{\sigma}[x(s)]\frac{\Delta}{\Delta x(s-1/2)}\left[  \frac{\nabla
y(s)}{\nabla x(s)}\right]  +\frac{1}{2}\widetilde{\tau}[x(s)]\left[
\frac{\Delta y(s)}{\Delta x(s)}+\frac{\nabla y(s)}{\nabla x(s)}\right]
+\lambda y(s)=0, \label{NUSeq}%
\end{equation}
where $\widetilde{\sigma}(x)$ and $\widetilde{\tau}(x)$ are polynomials of
degrees at most two and one in $x(s),$ respectively, $\lambda$ is a constant,
$\Delta y(s)=y(s+1)-y(s),$ $\nabla y(s)=y(s)-y(s-1),$ and $x(s)$ is a lattice
function\ that satisfies
\begin{equation}
\frac{x(s+1)+x(s)}{2}=\alpha x(s+\frac{1}{2})+\beta,\text{\hspace{0.1in}%
}\alpha,\beta\text{ are constants,} \label{condition1}%
\end{equation}%
\begin{equation}
x^{2}(s+1)+x^{2}(s)\text{ is a polynomial of degree at most two w.r.t.
}x(s+\frac{1}{2}). \label{condition2}%
\end{equation}
It should be pointed out that Eq. (\ref{NUSeq}) was obtained as a result of
approximating Eq. (\ref{hytype}) on a non-uniform lattice $x(s)$. Here, two
kinds of lattice functions $x(s)$ called \emph{non-uniform lattices} which
satisfy the conditions in Eqs. (\ref{condition1}) and (\ref{condition2}) are
\begin{equation}
x(s)=c_{1}q^{s}+c_{2}q^{-s}+c_{3}, \label{class1}%
\end{equation}%
\begin{equation}
x(s)=\widetilde{c}_{1}s^{2}+\widetilde{c}_{2}s+\widetilde{c}_{3},
\label{class2}%
\end{equation}
where $c_{i},\widetilde{c}_{i}$ are arbitrary constants and $c_{1}c_{2}\neq0$,
$\widetilde{c}_{1}\widetilde{c}_{2}\neq0$.

It was found that Eq. (\ref{NUSeq}) is of independent importance itself and
the equation intrigues many interesting questions. Its solutions essentially
generalize the solutions of the original differential equation and are of
interest in their own selves. Some of its solutions have been used in quantum
mechanics, the theory of group representations, and computational mathematics
\cite{nikiforov1991, nikiforov1988}. In particular, Suslov \cite{suslov1989}
established an analogous fundamental result for the difference equation on
non-uniform lattices, which generalizes the Rodrigues formula for polynomial
solutions on non-uniform lattices.

We should mention that the adjoint difference equation of Eq. (\ref{NUSeq})
for the case of non-uniform lattices is also of independent importance itself
and the adjoint equation may intrigue some other interesting questions. For
example, it could help us to obtain the particular solutions for the
difference equation of hypergeometric type in Eq. (\ref{NUSeq}) on non-uniform
lattices, or obtain an extension of the Rodrigues formula in the non-uniform
lattice case, etc. To our best knowledge, the adjoint difference equation of
Eq. (\ref{NUSeq}) for the case of uniform lattices such as $x(s)=s$ and
$x(s)=q^{s}$ has already been obtained in \cite{alvarez2011, area2003,
area2005}. However, for the case of non-uniform lattices where $x(s)$ is
defined in Eq. (\ref{class1}) or Eq. (\ref{class2}), it is more complicated
and difficult to establish and simplify and then solve the adjoint equations
on the non-uniform lattice,\ and as a result, the related study has not been
done yet.

Following the work of literature \cite{suslov1989}, the purpose of this
article is to establish an adjoint difference equation for the difference
equation in Eq. (\ref{NUSeq}) on non-uniform lattices where $x(s)$ is given in
Eq. (\ref{class1}) or Eq. (\ref{class2}). We will prove that the adjoint
difference equation is still a difference equation of hypergeometric type on
non-uniform lattices, and then obtain the particular solutions of the new
adjoint equation. In addition, we will prove another kind of fundamental
theorems for the Nikiforov-Uvarov-Suslov difference equation of hypergeometric
type in Eq. (\ref{NUSeq}), which are essentially new results and their
expressions are different from the Suslov theorem (as seen in the comparison
between Theorems 5.1-5.2 and Corollaries 5.1-5.2 in this article). As an
application of the new fundamental theorems, we use them to obtain the form of
particular solutions of the original difference equation of hypergeometric
type in Eq. (\ref{NUSeq}) on non-uniform lattices.

The rest of the paper is organized as follows. In section 2, we introduce some
preliminary information about the difference equation of hypergeometric type
on non-uniform lattices and give some related propositions and lemmas. In
section 3, we first consider a more general equation than Eq. (\ref{NUSeq}) on
non-uniform lattices and derive its adjoint equation. We then prove that the
adjoint equation is also a difference equation of hypergeometric type on
non-uniform lattices. As a special case of these results, we obtain the
adjoint equation of hypergeometric type for Eq. (\ref{NUSeq}) on non-uniform
lattices. In section 4, we derive the forms of particular solutions for the
general adjoint equation and its special adjoint equation for Eq.
(\ref{NUSeq}), respectively. In section 5, we use these particular solutions
to obtain particular solutions for both Eq. (\ref{NUSeq}) and its more general
difference equation of hypergeometric type on non-uniform lattices. In
addition, we give another two new fundamental theorems for both Eq.
(\ref{NUSeq}) and its more general difference equation of hypergeometric type
on non-uniform lattices. Finally, we give an example to illustrate the
application of the new fundamental theorems. The conclusion is then given in
section 6.

Throughout this paper, we follow those notations used in \cite{nikiforov1991,
nikiforov1988, suslov1989} which now have become the standard for analysis.
Furthermore, the properties listed below will be used in our study:
\begin{equation}
\frac{{x(z+\mu)+x(z)}}{2}=\alpha(\mu)x(z+\frac{\mu}{2})+\beta(\mu),
\label{condition3}%
\end{equation}%
\begin{equation}
x(z+\mu)-x(z)=\gamma(\mu)\nabla x(z+\frac{{\mu+1}}{2}), \label{condition4}%
\end{equation}
where
\begin{equation}
\alpha(\mu)=\left\{
{\begin{array}{*{20}c} {\frac{{q^{\frac{\mu } {2}} + q^{ - \frac{\mu } {2}} }} {2},} \\ {1;} \\ \end{array}}%
\right.  \beta(\mu)=\left\{
{\begin{array}{*{20}c} {\beta \frac{{1 - \alpha (\mu )}} {{1 - \alpha }},} \\ {\beta \mu ^2 ;} \\ \end{array}}%
\right.  \gamma(\mu)=\left\{
{\begin{array}{*{20}c} {\frac{{q^{\frac{\mu } {2}} - q^{ - \frac{\mu } {2}} }} {{q^{\frac{1} {2}} - q^{ - \frac{1} {2}} }},} \\ {\mu .} \\ \end{array}}%
\right.  \label{abc}%
\end{equation}

\section{{\protect\Large Preliminary Information and Lemmas}}

In this section, we give some preliminary information on the difference
equation of hypergeometric type on non-uniform lattices. Let $x(s)$ be a
lattice, where $s\in\mathbb{C}$ (complex numbers). It can be seen that for any
real $\nu$, $x_{\nu}(s)=x(s+\frac{\nu}{2})$ is also a lattice. Given a
function $f(s)$, we define two difference operators with respect to $x_{\nu
}(s)$ as%
\begin{equation}
\Delta_{\nu}f(s)=\frac{\Delta f(s)}{\Delta x_{\nu}(s)},\text{\hspace{0.1in}%
}\nabla_{\nu}f(s)=\frac{\nabla f(s)}{\nabla x_{\nu}(s)}.
\end{equation}
Moreover, for any nonnegative integer $n$, we define%
\[
\Delta_{\nu}^{(n)}f(s)=\left\{
\begin{aligned} &f(s), && \text{ if $n=0$,} \\ &\frac{\Delta}{\Delta x_{\nu+n-1}(s)}\cdots\frac{\Delta}{\Delta x_{{\nu}+1}(s)}\frac{\Delta}{\Delta x_{\nu}(s)}f(s), && ~\text{if $n\geq1$}. \end{aligned}\right.
\]%
\[
\nabla_{\nu}^{(n)}f(s)=\left\{
\begin{aligned} &f(s), && \text{ if $n=0$,} \\ &\frac{\nabla}{\nabla x_{\nu-n+1}}\cdots\frac{\nabla}{\nabla x_{\nu-1}(s)}\frac{\nabla}{\nabla x_{\nu}(s)}f(s), && ~\text{if $n\geq1$}. \end{aligned}\right.
\]

The following proposition can be verified straightforwardly.

\textbf{Proposition 2.1.} \textit{Given two functions} $f(s),g(s)$
\textit{with complex variable} $s$,\textit{\ the following difference
equalities hold}%
\begin{align*}
\Delta_{\nu}(f(s)g(s))  &  =f(s+1)\Delta_{\nu}g(s)+g(s)\Delta_{\nu}f(s)\\
&  =g(s+1)\Delta_{\nu}f(s)+f(s)\Delta_{\nu}g(s),
\end{align*}%
\begin{align*}
\Delta_{\nu}\left(  \frac{f(s)}{g(s)}\right)   &  =\frac{g(s+1)\Delta_{\nu
}f(s)-f(s+1)\Delta_{\nu}g(s)}{g(s)g(s+1)}\\
&  =\frac{g(s)\Delta_{\nu}f(s)-f(s)\Delta_{\nu}g(s)}{g(s)g(s+1)},
\end{align*}%
\begin{align*}
\nabla_{\nu}(f(s)g(s))  &  =f(s-1)\nabla_{\nu}g(s)+g(s)\nabla_{\nu}f(s)\\
&  =g(s-1)\nabla_{\nu}f(s)+f(s)\nabla_{\nu}g(s),
\end{align*}%
\begin{align*}
\nabla_{\nu}\left(  \frac{f(s)}{g(s)}\right)   &  =\frac{g(s-1)\nabla_{\nu
}f(s)-f(s-1)\nabla_{\nu}g(s)}{g(s)g(s-1)}\\
&  =\frac{g(s)\nabla_{\nu}f(s)-f(s)\nabla_{\nu}g(s)}{g(s)g(s-1)}.
\end{align*}

It can be seen that the difference equation of hypergeometric type in Eq.
(\ref{NUSeq}) can be written as
\begin{equation}
\widetilde{\sigma}[x(s)]\Delta_{-1}\nabla_{0}y(s)+\frac{\widetilde{\tau
}[x(s)]}{2}\left[  \Delta_{0}y(s)+\nabla_{0}y(s)\right]  +\lambda y(s)=0.
\label{NUSeq2}%
\end{equation}
Let
\[
z_{k}(s)=\Delta_{0}^{(k)}y(s)=\Delta_{k-1}\Delta_{k-2}\cdots\Delta
_{0}y(s),\quad k=1,2,....
\]
Then, for any nonnegative integer $k$, $z_{k}(s)$ satisfies an equation that
has the same type equation as Eq. (\ref{NUSeq2}) as \cite{nikiforov1991,
nikiforov1988}:
\begin{equation}
\widetilde{\sigma}_{k}[x_{k}(s)]\Delta_{k-1}\nabla_{k}z_{k}(s)+\frac
{\widetilde{\tau}_{k}[x_{k}(s)]}{2}\left[  \Delta_{k}z_{k}(s)+\nabla_{k}%
z_{k}(s)\right]  +\mu_{k}z_{k}(s)=0, \label{NUSeq3}%
\end{equation}
where $\mu_{k}$ is a constant, and $\widetilde{\sigma}_{k}(x_{k})$ and
$\widetilde{\tau}_{k}(x_{k})$ are polynomials of degrees at most two and one
in $x_{k}$, respectively, which are given as%
\begin{align*}
\widetilde{\sigma}_{k}[x_{k}(s)]  &  =\frac{\widetilde{\sigma}_{k-1}%
[x_{k-1}(s+1)]+\widetilde{\sigma}_{k-1}[x_{k-1}(s)]}{2}\\
&  +\frac{1}{4}\Delta_{k-1}\widetilde{\tau}_{k-1}(s)\frac{\Delta
x_{k}(s)+\nabla x_{k}(s)}{2\Delta x_{k-1}(s)}[\Delta x_{k-1}(s)]^{2}\\
&  +\frac{\widetilde{\tau}_{k-1}[x_{k-1}(s+1)]+\widetilde{\tau}_{k-1}%
[x_{k-1}(s)]}{2}\frac{\Delta x_{k}(s)-\nabla x_{k}(s)}{4},\\
\widetilde{\sigma}_{0}[x_{0}(s)]  &  =\widetilde{\sigma}[x(s)];
\end{align*}%
\begin{align*}
\widetilde{\tau}_{k}[x_{k}(s)]  &  =\Delta_{k-1}\widetilde{\sigma}%
_{k-1}\left[  x_{k-1}(s)\right]  +\Delta_{k-1}\widetilde{\tau}_{k-1}\left[
x_{k-1}(s)\right]  \frac{\Delta x_{k}(s)-\nabla x_{k}(s)}{4}\\
&  +\frac{\widetilde{\tau}_{k-1}[x_{k-1}(s+1)]+\widetilde{\tau}_{k-1}%
[x_{k-1}(s)]}{2}\frac{\Delta x_{k}(s)+\nabla x_{k}(s)}{2\Delta x_{k-1}(s)},\\
\widetilde{\tau}_{0}[x_{0}(s)]  &  =\widetilde{\tau}[x(s)];
\end{align*}%
\[
\mu_{k}=\mu_{k-1}+\Delta_{k-1}\widetilde{\tau}_{k-1}\left[  x_{k-1}(s)\right]
,~\mu_{0}=\lambda.
\]

To analyze additional properties of solutions of Eq. (\ref{NUSeq3}), it is
convenient to use the equality
\[
\frac{1}{2}\left[  \Delta_{k}z_{k}(s)+\nabla_{k}z_{k}(s)\right]  =\Delta
_{k}z_{k}(s)-\frac{1}{2}\Delta\left[  \nabla_{k}z_{k}(s)\right]
\]
and rewrite Eq. (\ref{NUSeq3}) in an equivalent expression as
\begin{equation}
\sigma_{k}(s)\Delta_{k-1}\nabla_{k}z_{k}(s)+\tau_{k}(s)\Delta_{k}z_{k}%
(s)+\mu_{k}z_{k}(s)=0,
\end{equation}
where%
\begin{align}
&  \sigma_{k}(s)=\widetilde{\sigma}_{k}[x_{k}(s)]-\frac{1}{2}\widetilde{\tau
}_{k}[x_{k}(s)]\nabla x_{k+1}(s),\label{19}\\
&  \sigma_{k}(s)=\sigma_{k-1}(s)=\sigma(s);\label{20}\\
\tau_{k}(s)\nabla x_{k+1}(s)  &  =\sigma(s+k)-\sigma(s)+\tau(s+k)\nabla
x_{1}(s+k),\label{21}\\
\tau_{k}(s)  &  =\widetilde{\tau}_{k}[x_{k}(s)];\\
\mu_{k}  &  =\lambda+\sum_{j=0}^{k-1}\Delta_{j}\tau_{j}(s).
\end{align}
Eq. (\ref{NUSeq3}) can be further rewritten into a self-adjoint form as
\[
\Delta_{k-1}\left[  \sigma(s)\rho_{k}(s)\nabla_{k}z_{k}(s)\right]  +\mu
_{k}\rho_{k}(s)z_{k}(s)=0,
\]
where $\rho_{k}(s)$ satisfies the Pearson type difference equation as
\[
\Delta_{k-1}\left[  \sigma(s)\rho_{k}(s)\right]  =\tau_{k}(s)\rho_{k}(s).
\]
Letting $\rho(s)=\rho_{0}(s)$, then we have
\[
\rho_{k}(s)=\rho(s+k)\prod_{i=1}^{k}\sigma(s+i).
\]
Thus, if for a positive integer $n$ such that $\mu_{n}=0$ and
\begin{equation}
\lambda=\lambda_{n}:=-\sum_{j=0}^{n-1}\Delta_{j}\tau_{j}(s),~\text{and}%
~\lambda_{m}\neq\lambda_{n}~\text{for}~m=0,1,\dots,n-1,
\end{equation}
then Eq. (\ref{NUSeq3}) has a polynomial solution $y_{n}[x(s)]$ of degree $n$
about $x(s)$, which can be expressed by the difference analog of the Rodrigues
formula \cite{nikiforov1991, nikiforov1988}:
\[
y_{n}[x(s)]=\frac{1}{\rho(s)}\nabla_{n}^{(n)}\left[  \rho_{n}(s)\right]
=\frac{1}{\rho(s)}\Delta_{-n}^{(n)}\left[  \rho_{n}(s-n)\right]  .
\]

Furthermore, when $k$ is nonnegative integer, S.K. Suslov in \cite{suslov1989}
gave the following extension definitions of equalities (\ref{19}) and
(\ref{21}).

\textbf{Definition 2.1. }\textit{Let} $\nu\in R,x=x(z)$ \textit{be a lattice
satisfying the two conditions given in Eqs. (\ref{condition3}%
)-(\ref{condition4})}. \textit{Then, functions} $\tilde{\sigma}_{\nu}(z)$
\textit{and} $\tau_{\nu}(z)$ \textit{are defined by the equalities}
\begin{equation}
\tilde{\sigma}_{\nu}(z)=\sigma(z)+\frac{1}{2}\tau_{\nu}(z)\nabla x_{\nu+1}(z),
\label{22}%
\end{equation}%
\begin{equation}
\tau_{\nu}(z)\nabla x_{\nu+1}(z)=\sigma(z+\nu)-\sigma(z)+\tau(z+\nu)\nabla
x_{1}(z+\nu). \label{23}%
\end{equation}
S.K. Suslov further studied the following extension of the
Nikiforov-Uvarov-Suslov equation in Eq. (\ref{NUSeq}) as:%

\begin{equation}
\sigma(z)\frac{\Delta}{{\Delta x_{\nu-\mu-1}(z)}}(\frac{{\nabla y(z)}}{{\nabla
x_{\nu-\mu}(z)}})+\tau_{\nu-\mu}(z)\frac{{\Delta y(z)}}{{\Delta x_{\nu-\mu
}(z)}}+\lambda y(z)=0 \label{NUSeq1}%
\end{equation}
where $\nu,\mu\in R,x=x(z)$ be a lattice satisfying the two conditions given
in Eqs. (\ref{condition3})-(\ref{condition4}), and obtained several important
results, which can be seen in the books \cite{nikiforov1991, nikiforov1988}.

\textbf{Lemma 2.1 }\cite{suslov1989}\textbf{.} $\tilde{\sigma}_{\nu}%
(z)=\tilde{\sigma}_{\nu}(x_{\nu})$ \textit{and} $\tau_{\nu}(z)=\tau_{\nu
}(x_{\nu})$ \textit{are polynomials of degrees at most two and one,
respectively, in the variable} $x_{\nu}(s)=x(s+\frac{\nu}{2}).$

\textbf{Lemma 2.2 }\cite{suslov1989}\textbf{.} \textit{Under the hypotheses of
Lemma 2.1, the function}%
\[
Q(s)=\nu(\mu)\sigma(s)-\tau_{\nu}(s)[x_{\nu-\mu}(s)-x_{\nu-\mu}(z)]
\]
\textit{has the form of}%
\[
Q(s)=A+B[x_{\nu}(s)-x_{\nu}(z)]+C[x_{\nu}(s)-x_{\nu}(z)][x_{\nu}(s)-x_{\nu
}(z-\mu)],
\]
\textit{where}%
\[
A=\nu(\mu)\sigma(z),\hspace{0.1in}B=-\tau_{\nu-\mu}(z),\hspace{0.1in}%
C=-\kappa_{\mu-2\nu+1}.
\]

The following identities about the explicit form of $\tau_{\nu}(s),\mu_{k}$
and{\Large \ }$\lambda_{n}$ are not difficult to check when the non-uniform
lattice is either $x(s)=c_{1}q^{s}+c_{2}q^{-s}+c_{3}$ or $x(s)=\tilde{c}%
_{1}s^{2}+\tilde{c}_{2}s+\tilde{c}_{3}$.

\textbf{Proposition 2.2} \cite{suslov1989}\textbf{.} \textit{Given any real}
$\nu,$ let $\alpha(\nu)$ \textit{and} $\nu(\nu)$ \textit{be defined in Eq.
(\ref{abc}),} $\kappa_{\nu}$ \textit{be defined in Eq. (\ref{tauka})},
\textit{and if} $x(s)=c_{1}q^{s}+c_{2}q^{-s}+c_{3}$, \textit{then}
\[
\tau_{\nu}(s)=\kappa_{2\nu+1}x_{\nu}(s)+c(\nu),
\]
\textit{where} $c(\nu)$ \textit{is a function with respect to} $\nu$
\textit{as}%
\begin{align*}
c(\nu)  &  =c_{3}(1-q^{\frac{\nu}{2}})(q^{\frac{\nu}{2}}-q^{-\nu}%
)+c_{3}(2-q^{\frac{\nu}{2}}-q^{-\frac{\nu}{2}})\nu(\nu)+2\widetilde{\tau
}(0)\alpha(\nu)\\
&  +\widetilde{\sigma}(0)\nu(\nu).
\end{align*}

\textit{On the other hand, if} $x(s)=\widetilde{c}_{1}s^{2}+\widetilde{c}%
_{2}s+\widetilde{c}_{3}$,\textit{\ then}
\begin{equation}
\tau_{\nu}(s)=\kappa_{2\nu+1}x_{\nu}(s)+\widetilde{c}(\nu),\nonumber
\end{equation}
\textit{where }$\widetilde{c}(\nu)$ \textit{is a function with respect to}
$\nu$ \textit{as}%
\[
\widetilde{c}(\nu)=\frac{\widetilde{\sigma}"}{4}\widetilde{c}_{1}{\nu}%
^{3}+\frac{3\widetilde{\tau}^{\prime}}{4}\widetilde{c}_{1}{\nu}^{2}%
+\widetilde{\sigma}(0)\nu+2\widetilde{\tau}(0).
\]

\textbf{Proposition 2.3} \cite{suslov1989}\textbf{.} \textit{For} $\alpha
(\mu)$ \textit{and} $\nu(\mu)$, \textit{it holds that}
\[
\sum\limits_{j=0}^{k-1}{\alpha(2j)}=\alpha(k-1)\nu(k),\text{ \hspace{0.1in}%
}\sum\limits_{j=0}^{k-1}{\nu(2j)}=\nu(k-1)\nu(k).
\]

\textbf{Proposition 2.4} \cite{nikiforov1991, nikiforov1988}\textbf{.}
\textit{For} $x(s)=c_{1}q^{s}+c_{2}q^{-s}+c_{3}$ \textit{or}
$x(s)=\widetilde{c}_{1}s^{2}+\widetilde{c}_{2}s+\widetilde{c}_{3}$, \textit{it
holds that}%
\begin{equation}
\mu_{k}=\lambda+\kappa_{k}\nu(k),\quad k=1,2,...,n,
\end{equation}%
\begin{equation}
\lambda_{n}=-n\kappa_{n}.
\end{equation}

\section{{\protect\Large Adjoint difference equation}}

We now seek the second-order adjoint difference equation corresponding to Eq.
(\ref{NUSeq}). To this end, we first consider the operator%
\begin{equation}
L[y(z)]\equiv\sigma(z)\frac{\Delta}{{\Delta x_{\nu-\mu-1}(z)}}(\frac{{\nabla
y(z)}}{{\nabla x_{\nu-\mu}(z)}})+\tau_{\nu-\mu}(z)\frac{{\Delta y(z)}}{{\Delta
x_{\nu-\mu}(z)}}+\lambda y(z). \label{operL}%
\end{equation}
One may see that the equation%
\begin{equation}
L[y(z)]\equiv\sigma(z)\frac{\Delta}{{\Delta x_{\nu-\mu-1}(z)}}(\frac{{\nabla
y(z)}}{{\nabla x_{\nu-\mu}(z)}})+\tau_{\nu-\mu}(z)\frac{{\Delta y(z)}}{{\Delta
x_{\nu-\mu}(z)}}+\lambda y(z)=0 \label{NUSeq4}%
\end{equation}
is a more generalized Nikiforov-Uvarov-Suslov equation than Eq. (\ref{NUSeq}),
and it can be reduced to Eq. (\ref{NUSeq}) by letting $\mu=\nu.$ If we rewrite
Eq. (\ref{NUSeq4}) as
\begin{equation}
\tilde{\sigma}_{\nu-\mu}(z)\frac{\Delta}{{\Delta x_{\nu-\mu-1}(z)}}%
(\frac{{\nabla y(z)}}{{\nabla x_{\nu-\mu}(z)}})+\frac{{\tau_{\nu-\mu}(z)}}%
{2}[\frac{{\Delta y(z)}}{{\Delta x_{\nu-\mu}(z)}}+\frac{{\nabla y(z)}}{{\nabla
x_{\nu-\mu}(z)}}]+\lambda y(z)=0, \label{NUSeq5}%
\end{equation}
where
\[
\tilde{\sigma}_{\nu-\mu}(z)=\sigma(z)+\frac{1}{2}\tau_{\nu-\mu}(z)\nabla
x_{\nu-\mu+1}(z),
\]
by Lemma 2.1, one may see that $\tilde{\sigma}_{\nu-\mu}(z)$ and $\tau
_{\nu-\mu}(z)$ are polynomials of degrees at most two and one, respectively,
in the variable $x_{\nu-\mu}(s)$, and hence, Eq. (\ref{NUSeq5}) is a
difference equation of hypergeometric type.

\textbf{Definition 3.1. }\textit{For }$y(z)$\textit{ and }$w(z)$\textit{, the
scalar product }$\langle w(z),y(z)\rangle$\textit{ with respect to }$\Delta
x_{\nu-\mu-1}(z)$\textit{ is defined as}%
\[
\langle w(z),y(z)\rangle=\sum_{z=a}^{b-1}w(z)y(z){\Delta x_{\nu-\mu-1}(z),}%
\]
\textit{where }$a,b$\textit{ are complex with the same imaginary parts, and
}$b-a\in N.$

\textbf{Definition 3.2. }\textit{For }$w(z)$\textit{ and operator }%
$L[y(z)]$\textit{, assume that the boundary conditions }%
$w(a)=w(b)=0,y(a)=y(b)=0$\textit{ are satisfied. If the scalar product }%
\[
\langle w(z),L[y(z)]\rangle=\langle y(z),L^{\ast}[w(z)]\rangle
\]
\textit{holds, then the operator }$L^{\ast}[w(z)]$\textit{ is called the
adjoint operator of }$L[y(z)]$\textit{, and }$L^{\ast}[w(z)]=0$\textit{ is
called the adjoint equation of }$L[y(z)]=0.$

We now find the operator $L^{\ast}[w(z)]$. Since%
\begin{align*}
&  \langle w(z),L[y(z)]\rangle\\
&  =\sum_{z=a}^{b-1}w(z)L[y(z)]{\Delta x_{\nu-\mu-1}(z)}\\
&  =\sum_{z=a}^{b-1}w(z)\{\sigma(z)\frac{\Delta}{{\Delta x_{\nu-\mu-1}(z)}%
}(\frac{{\nabla y(z)}}{{\nabla x_{\nu-\mu}(z)}})+\tau_{\nu-\mu}(z)\frac
{{\Delta y(z)}}{{\Delta x_{\nu-\mu}(z)}}+\lambda y(z)\}{\Delta x_{\nu-\mu
-1}(z),}%
\end{align*}
using the summation by parts and the boundary conditions, we obtain%
\begin{align*}
&  \sum_{z=a}^{b-1}w(z)\sigma(z)\frac{\Delta}{{\Delta x_{\nu-\mu-1}(z)}}%
(\frac{{\nabla y(z)}}{{\nabla x_{\nu-\mu}(z)}}){\Delta x_{\nu-\mu-1}(z)}\\
&  =\sum_{z=a}^{b-1}w(z)\sigma(z)\Delta(\frac{{\nabla y(z)}}{{\nabla
x_{\nu-\mu}(z)}})=-\sum_{z=a}^{b-1}\frac{{\nabla y(z)}}{{\nabla x_{\nu-\mu
}(z)}}\nabla(w(z)\sigma(z))\\
&  =-\sum_{z=a}^{b-1}{\nabla y(z)}\frac{\nabla(w(z)\sigma(z))}{{\nabla
x_{\nu-\mu}(z)}}=\sum_{z=a}^{b-1}{y(z)\Delta\lbrack}\frac{\nabla
(w(z)\sigma(z))}{{\nabla x_{\nu-\mu}(z)}}]\\
&  =\sum_{z=a}^{b-1}{y(z)}\frac{\Delta}{{\Delta x_{\nu-\mu-1}(z)}}{[}%
\frac{\nabla(w(z)\sigma(z))}{{\nabla x_{\nu-\mu}(z)}}]{\Delta x_{\nu-\mu
-1}(z),}%
\end{align*}
and%
\begin{align*}
&  \sum_{z=a}^{b-1}w(z)\tau_{\nu-\mu}(z)\frac{{\Delta y(z)}}{{\Delta
x_{\nu-\mu}(z)}}{\Delta x_{\nu-\mu-1}(z)}\\
&  =\sum_{z=a}^{b-1}w(z)\tau_{\nu-\mu}(z)\frac{{\Delta x_{\nu-\mu-1}(z)}%
}{{\Delta x_{\nu-\mu}(z)}}{\Delta y(z)}=-\sum_{z=a}^{b-1}{y(z)\nabla\lbrack
}w(z)\tau_{\nu-\mu}(z)\frac{{\Delta x_{\nu-\mu-1}(z)}}{{\Delta x_{\nu-\mu}%
(z)}}]\\
&  =-\sum_{z=a}^{b-1}{y(z)}\frac{{\nabla}}{{\Delta x_{\nu-\mu-1}(z)}}%
{[}w(z)\tau_{\nu-\mu}(z)\frac{{\Delta x_{\nu-\mu-1}(z)}}{{\Delta x_{\nu-\mu
}(z)}}]{\Delta x_{\nu-\mu-1}(z).}%
\end{align*}
Thus, we let%
\begin{align*}
&  \langle w(z),L[y(z)]\rangle\\
&  =\sum_{z=a}^{b-1}y(z)\{\frac{\Delta}{{\Delta x_{\nu-\mu-1}(z)}}{[}%
\frac{\nabla(w(z)\sigma(z))}{{\nabla x_{\nu-\mu}(z)}}]-\frac{{\nabla}}{{\Delta
x_{\nu-\mu-1}(z)}}{[}w(z)\tau_{\nu-\mu}(z)\frac{{\Delta x_{\nu-\mu-1}(z)}%
}{{\Delta x_{\nu-\mu}(z)}}]\\
&  -\lambda w(z)\}{\Delta x_{\nu-\mu-1}(z)}\\
&  =\sum_{z=a}^{b-1}y(z)L^{\ast}[w(z)]{\Delta x_{\nu-\mu-1}(z)}\\
&  {=}\langle y(z),L^{\ast}[w(z)]\rangle,
\end{align*}
which gives%
\begin{equation}
L^{\ast}[w(z)]=\frac{\Delta}{{\Delta x_{\nu-\mu-1}(z)}}{[}\frac{\nabla
(w(z)\sigma(z))}{{\nabla x_{\nu-\mu}(z)}}]-\frac{{\nabla}}{{\Delta x_{\nu
-\mu-1}(z)}}{[}w(z)\tau_{\nu-\mu}(z)\frac{{\Delta x_{\nu-\mu-1}(z)}}{{\Delta
x_{\nu-\mu}(z)}}]-\lambda w(z). \label{operLstar}%
\end{equation}
Therefore, we define Eq. (\ref{operLstar}) as the \textbf{adjoint operator }of
Eq. (\ref{operL})\textbf{. }

It can be seen that%
\begin{align}
&  \frac{\Delta}{{\Delta x_{\nu-\mu-1}(z)}}{[}\frac{\nabla(w(z)\sigma
(z))}{{\nabla x_{\nu-\mu}(z)}}]\nonumber\\
&  =\frac{\Delta}{{\Delta x_{\nu-\mu-1}(z)}}[\sigma(z-1)\frac{\nabla
w(z)}{{\nabla x_{\nu-\mu}(z)}}+w(z)\frac{\nabla\sigma(z)}{{\nabla x_{\nu-\mu
}(z)}}\nonumber\\
&  =\sigma(z-1)\frac{\Delta}{{\Delta x_{\nu-\mu-1}(z)}}[\frac{\nabla
w(z)}{{\nabla x_{\nu-\mu}(z)}}]+\frac{\Delta w(z)}{{\Delta x_{\nu-\mu}(z)}%
}\frac{\Delta\sigma(z-1)}{{\Delta x_{\nu-\mu}(z)}}\nonumber\\
&  +\frac{\Delta\sigma(z)}{{\Delta x_{\nu-\mu}(z)}}\frac{\Delta w(z)}{{\Delta
x_{\nu-\mu-1}(z)}}+w(z)\frac{\Delta}{{\Delta x_{\nu-\mu-1}(z)}}[\frac
{\nabla\sigma(z)}{{\nabla x_{\nu-\mu}(z)}}], \label{31}%
\end{align}%
\begin{align}
&  \frac{{\nabla}}{{\Delta x_{\nu-\mu-1}(z)}}{[}w(z)\tau_{\nu-\mu}%
(z)\frac{{\Delta x_{\nu-\mu-1}(z)}}{{\Delta x_{\nu-\mu}(z)}}]\nonumber\\
&  =\tau_{\nu-\mu}(z-1)\frac{{\Delta x_{\nu-\mu-1}(z)}}{{\Delta x_{\nu-\mu
}(z)}}\frac{{\nabla w(z)}}{{\Delta x_{\nu-\mu-1}(z)}}+{w(z)}\frac{{\nabla}%
}{{\Delta x_{\nu-\mu-1}(z)}}{[}\tau_{\nu-\mu}(z)\frac{{\Delta x_{\nu-\mu
-1}(z)}}{{\Delta x_{\nu-\mu}(z)}}], \label{32}%
\end{align}%
\begin{align}
\frac{{\nabla w(z)}}{{\nabla x_{\nu-\mu}(z)}}  &  =\frac{{\Delta w(z)}%
}{{\Delta x_{\nu-\mu}(z)}}-\Delta\lbrack\frac{{\nabla w(z)}}{{\nabla
x_{\nu-\mu}(z)}}]\nonumber\\
&  =\frac{{\Delta w(z)}}{{\Delta x_{\nu-\mu}(z)}}-\frac{\Delta}{{\Delta
x_{\nu-\mu-1}(z)}}[\frac{{\nabla w(z)}}{{\nabla x_{\nu-\mu}(z)}}]{\Delta
x_{\nu-\mu-1}(z).} \label{33}%
\end{align}
Substituting Eqs. (\ref{31}-\ref{33}) into Eq. (\ref{operLstar}), we obtain
another expression of $L^{\ast}[w(z)]$ as
\begin{equation}
L^{\ast}[w(z)]\equiv\sigma^{\ast}(z)\frac{\Delta}{{\Delta x_{\nu-\mu-1}(z)}%
}(\frac{{\nabla w(z)}}{{\nabla x_{\nu-\mu}(z)}})+\tau_{\nu-\mu}^{\ast}%
(z)\frac{{\Delta w(z)}}{{\Delta x_{\nu-\mu}(z)}}+\lambda_{\nu-\mu}^{\ast}w(z),
\label{operLstar2}%
\end{equation}
where
\begin{equation}
\sigma^{\ast}(z)=\sigma(z-1)+\tau_{\nu-\mu}(z-1)\nabla x_{\nu-\mu-1}(z),
\label{41a1}%
\end{equation}%
\begin{equation}
{\tau_{\nu-\mu}^{\ast}(z)=\frac{{\sigma(z+1)-\sigma(z-1)-\tau_{\nu-\mu
}(z-1)\nabla x_{\nu-\mu-1}(z)}}{{\Delta x_{\nu-\mu-1}(z)}},} \label{41b1}%
\end{equation}%
\begin{equation}
{\lambda_{\nu-\mu}^{\ast}=\lambda-\frac{\Delta}{{\Delta x_{\nu-\mu-1}(z)}%
}(\frac{{\tau_{\nu-\mu}(z-1)\nabla x_{\nu-\mu-1}(z)-\nabla\sigma(z)}}{{\nabla
x_{\nu-\mu}(z)}}).} \label{41c1}%
\end{equation}

We now seek the relationship between $\rho_{\nu-\mu}(z)L[y(z)]$ and $L^{\ast
}[\rho_{\nu-\mu}(z)y(z)]$. It can be seen that $\rho_{\nu-\mu}(z)L[y(z)]$ has
a self-adjoint form as
\begin{equation}
\rho_{\nu-\mu}(z)L[y(z)]=\frac{\Delta}{{\Delta x_{\nu-\mu-1}(z)}}%
(\sigma(z)\rho_{\nu-\mu}(z)\frac{{\nabla y(z)}}{{\nabla x_{\nu-\mu}(z)}%
})+\lambda\rho_{\nu-\mu}(z)y(z), \label{selfadj}%
\end{equation}
where $\rho_{\nu-\mu}(z)$ satisfies the Pearson type equation
\begin{equation}
\frac{{\Delta(\sigma(z)\rho_{\nu-\mu}(z))}}{{\Delta x_{\nu-\mu-1}(z)}}%
=\tau_{\nu-\mu}(z)\rho_{\nu-\mu}(z). \label{ptype}%
\end{equation}

Introducing $w(z)=\rho_{\nu-\mu}(z)y(z)$, we obtain
\[
\frac{{\nabla w(z)}}{{\nabla x_{\nu-\mu}(z)}}=\rho_{\nu-\mu}(z)\frac{{\nabla
y(z)}}{{\nabla x_{\nu-\mu}(z)}}+\frac{{\nabla\rho_{\nu-\mu}(z)}}{{\nabla
x_{\nu-\mu}(z)}}y(z-1),
\]
that is,
\begin{equation}
\rho_{\nu-\mu}(z)\frac{{\nabla y(z)}}{{\nabla x_{\nu-\mu}(z)}}=\frac{{\nabla
w(z)}}{{\nabla x_{\nu-\mu}(z)}}-\frac{{\nabla\rho_{\nu-\mu}(z)}}{{\nabla
x_{\nu-\mu}(z)}}y(z-1). \label{eq36}%
\end{equation}
From Eq. (\ref{ptype}), we have
\[
\Delta(\sigma(z)\rho_{\nu-\mu}(z))=\tau_{\nu-\mu}(z)\rho_{\nu-\mu}(z)\Delta
x_{\nu-\mu-1}(z),
\]
implying
\[
\rho_{\nu-\mu}(z)\Delta\sigma(z)+\Delta\rho_{\nu-\mu}(z)\sigma(z+1)=\tau
_{\nu-\mu}(z)\rho_{\nu-\mu}(z)\Delta x_{\nu-\mu-1}(z),
\]%
\[
\Delta\rho_{\nu-\mu}(z)=\frac{{\tau_{\nu-\mu}(z)\Delta x_{\nu-\mu-1}%
(z)-\Delta\sigma(z)}}{{\sigma(z+1)}}\rho_{\nu-\mu}(z),
\]
and hence
\begin{align}
\frac{{\nabla\rho_{\nu-\mu}(z)}}{{\nabla x_{\nu-\mu}(z)}}y(z-1)  &
=\frac{{\Delta\rho_{\nu-\mu}(z-1)y(z-1)}}{{\nabla x_{\nu-\mu}(z)}}\nonumber\\
&  =\frac{{\tau_{\nu-\mu}(z-1)\nabla x_{\nu-\mu-1}(z)-\nabla\sigma(z)}%
}{{\sigma(z)\nabla x_{\nu-\mu}(z)}}w(z-1). \label{eq37}%
\end{align}
Substituting Eqs. (\ref{eq36})-(\ref{eq37}) into Eq. (\ref{selfadj}), we
obtain%
\begin{align*}
&  \hspace{-0.2in}\rho_{\nu-\mu}(z)L[y(z)]=\frac{\Delta}{{\Delta x_{\nu-\mu
-1}(z)}}[\sigma(z)\frac{{\nabla w(z)}}{{\nabla x_{\nu-\mu}(z)}}]\\
&  -\frac{\Delta}{{\Delta x_{\nu-\mu-1}(z)}}[\frac{{\tau_{\nu-\mu}(z-1)\nabla
x_{\nu-\mu-1}(z)-\nabla\sigma(z)}}{{\nabla x_{\nu-\mu}(z)}}w(z-1)]+\lambda
w(z).
\end{align*}
This gives%
\begin{align*}
\hspace{-0.2in}\rho_{\nu-\mu}(z)L[y(z)]  &  =\sigma(z)\frac{\Delta}{{\Delta
x_{\nu-\mu-1}(z)}}[\frac{{\nabla w(z)}}{{\nabla x_{\nu-\mu}(z)}}%
]+\frac{{\Delta\sigma(z)}}{{\Delta x_{\nu-\mu-1}(z)}}\frac{{\Delta w(z)}%
}{{\Delta x_{\nu-\mu}(z)}}\hfill\\
&  -\frac{{\tau_{\nu-\mu}(z-1)\nabla x_{\nu-\mu-1}(z)-\nabla\sigma(z)}%
}{{\nabla x_{\nu-\mu}(z)}}\frac{{\nabla w(z)}}{{\Delta x_{\nu-\mu-1}(z)}%
}\hfill\\
&  -\frac{\Delta}{{\Delta x_{\nu-\mu-1}(z)}}[\frac{{\tau_{\nu-\mu}(z-1)\nabla
x_{\nu-\mu-1}(z)-\nabla\sigma(z)}}{{\nabla x_{\nu-\mu}(z)}}]w(z)+\lambda w(z),
\end{align*}
which implies that%
\begin{align}
\rho_{\nu-\mu}(z)L[y(z)]  &  =\sigma(z)\frac{\Delta}{{\Delta x_{\nu-\mu-1}%
(z)}}(\frac{{\nabla w(z)}}{{\nabla x_{\nu-\mu}(z)}})+\frac{{\Delta\sigma(z)}%
}{{\Delta x_{\nu-\mu-1}(z)}}\frac{{\Delta w(z)}}{{\Delta x_{\nu-\mu}(z)}%
}\nonumber\\
&  -\frac{{\tau_{\nu-\mu}(z-1)\nabla x_{\nu-\mu-1}(z)-\nabla\sigma(z)}%
}{{\Delta x_{\nu-\mu-1}(z)}}\frac{{\nabla w(z)}}{{\nabla x_{\nu-\mu}(z)}%
}\nonumber\\
&  -\frac{\Delta}{{\Delta x_{\nu-\mu-1}(z)}}[\frac{{\tau_{\nu-\mu}(z-1)\nabla
x_{\nu-\mu-1}(z)-\nabla\sigma(z)}}{{\nabla x_{\nu-\mu}(z)}}]w(z)+\lambda w(z).
\label{39}%
\end{align}
Using the following difference equalities%
\begin{align}
\frac{{\Delta w(z)}}{{\Delta x_{\nu-\mu}(z)}}-\frac{{\nabla w(z)}}{{\nabla
x_{\nu-\mu}(z)}}=  &  \Delta(\frac{{\nabla w(z)}}{{\nabla x_{\nu-\mu}(z)}%
}),\label{39a}\\
\frac{{\nabla w(z)}}{{\nabla x_{\nu-\mu}(z)}}=  &  \frac{{\Delta w(z)}%
}{{\Delta x_{\nu-\mu}(z)}}-\Delta(\frac{{\nabla w(z)}}{{\nabla x_{\nu-\mu}%
(z)}}), \label{39b}%
\end{align}
we can simplify Eq. (\ref{39}) to
\begin{align}
&  \rho_{\nu-\mu}(z)L[y(z)]\nonumber\\
&  =\sigma^{\ast}(z)\frac{\Delta}{{\Delta x_{\nu-\mu-1}(z)}}(\frac{{\nabla
w(z)}}{{\nabla x_{\nu-\mu}(z)}})+\tau_{\nu-\mu}^{\ast}(z)\frac{{\Delta w(z)}%
}{{\Delta x_{\nu-\mu}(z)}}+\lambda_{\nu-\mu}^{\ast}w(z), \label{41}%
\end{align}
where
\begin{equation}
\sigma^{\ast}(z)=\sigma(z-1)+\tau_{\nu-\mu}(z-1)\nabla x_{\nu-\mu-1}(z),
\label{41a}%
\end{equation}%
\begin{equation}
{\tau_{\nu-\mu}^{\ast}(z)=\frac{{\sigma(z+1)-\sigma(z-1)-\tau_{\nu-\mu
}(z-1)\nabla x_{\nu-\mu-1}(z)}}{{\Delta x_{\nu-\mu-1}(z)}},} \label{41b}%
\end{equation}%
\begin{equation}
{\lambda_{\nu-\mu}^{\ast}=\lambda-\frac{\Delta}{{\Delta x_{\nu-\mu-1}(z)}%
}(\frac{{\tau_{\nu-\mu}(z-1)\nabla x_{\nu-\mu-1}(z)-\nabla\sigma(z)}}{{\nabla
x_{\nu-\mu}(z)}}).} \label{41c}%
\end{equation}
Comparing with Eq. (\ref{operLstar2}), we see that the right-hand-side of Eq.
(\ref{41}) is $L^{\ast}[w(z)]$. Hence, we have obtained an important
relationship between the adjoint difference operator and the original
difference operator as described in Proposition 3.1.

\textbf{Proposition 3.1. }\textit{For }$y(z)$\textit{, it holds }%
\begin{equation}
L^{\ast}[\rho_{\nu-\mu}(z)y(z)]=\rho_{\nu-\mu}(z)L[y(z)]. \label{42}%
\end{equation}
Thus, we define
\begin{equation}
L^{\ast}[w(z)]=\sigma^{\ast}(z)\frac{\Delta}{{\Delta x_{\nu-\mu-1}(z)}}%
(\frac{{\nabla w(z)}}{{\nabla x_{\nu-\mu}(z)}})+\tau_{\nu-\mu}^{\ast}%
(z)\frac{{\Delta w(z)}}{{\Delta x_{\nu-\mu}(z)}}+\lambda_{\nu-\mu}^{\ast
}w(z)=0 \label{42a}%
\end{equation}
as the \textbf{adjoint difference equation}\emph{\ }corresponding to Eq.
(\ref{NUSeq4}). In particular, when $\mu=\nu$, Eq. (\ref{42a}) gives the
\textbf{adjoint difference equation}\emph{\ }corresponding to Eq.
(\ref{NUSeq}). Based on Definition 2.1 and Proposition 2.2, it is not
difficult to obtain the following corollary.

\textbf{Corollary 3.1.} \textit{Eqs. (\ref{41b}) and (\ref{41c}) can be
simplified to}
\begin{equation}
\tau_{\nu-\mu}^{\ast}(z)=-\tau_{\nu-\mu-2}(z+1)=-\kappa_{2\nu-2\mu-3}%
x_{\nu-\mu}(z)+c(\nu-\mu), \label{43a}%
\end{equation}%
\begin{equation}
\lambda_{\nu-\mu}^{\ast}=\lambda-\Delta_{\nu-\mu-1}\tau_{\nu-\mu-1}%
(z)=\lambda-\kappa_{2\nu-2\mu-1}. \label{43b}%
\end{equation}

\textbf{Proof.} Since%
\begin{align*}
&  \hspace{-0.2in}\sigma(z+1)-\sigma(z-1)-\tau_{\nu-\mu}(z-1)\nabla x_{\nu
-\mu-1}(z)\\
&  =\sigma(z+1)-\sigma(z-1)-\tau_{\nu-\mu}(z-1)\nabla x_{\nu-\mu+1}(z-1)\\
&  =\sigma(z+1)-\sigma(z-1+\nu-\mu)-\tau(z-1+\nu-\mu)\nabla x_{1}(z-1+\nu
-\mu)\\
&  =-\tau_{\nu-\mu-2}(z+1)\nabla x_{\nu-\mu-1}(z+1)\\
&  =-\tau_{\nu-\mu-2}(z+1)\Delta x_{\nu-\mu-1}(z),
\end{align*}
we obtain from Eq. (\ref{41b}) and Proposition 2.2 that%
\begin{align*}
\tau_{\nu-\mu}^{\ast}(z)  &  =-\tau_{\nu-\mu-2}(z+1)\\
&  =-\kappa_{2\nu-2\mu-3}x_{\nu-\mu-2}(z+1)+c(\nu-\mu)\\
&  =-\kappa_{2\nu-2\mu-3}x_{\nu-\mu}(z)+c(\nu-\mu).
\end{align*}
Using a similar argument, we have%
\[
\tau_{\nu-\mu}(z-1)\nabla x_{\nu-\mu-1}(z)-\nabla\sigma(z)=\tau_{\nu-\mu
}(z-1)\nabla x_{\nu-\mu-1}(z)+\sigma(z-1)-\sigma(z),
\]%
\[
\tau(z-1+\nu-\mu)\nabla x_{1}(z-1+\nu-\mu)+\sigma(z-1+\nu-\mu)-\sigma
(z)=\tau_{\nu-\mu-1}(z)\nabla x_{\nu-\mu}(z).
\]
Hence, we obtain%
\begin{align*}
\lambda_{\nu-\mu}^{\ast}  &  =\lambda-\Delta_{\nu-\mu-1}\tau_{\nu-\mu-1}(z)\\
&  =\lambda-\Delta_{\nu-\mu-1}\{\kappa_{2\nu-\mu-1}x_{\nu-\mu-1}(z)\}\\
&  =\lambda-\kappa_{2\nu-2\mu-1},
\end{align*}
and complete the proof.

Regarding to the adjoint difference equation in Eq. (\textit{\ref{42a}}), we
find some interesting dual properties as described in the following proposition.

\textbf{Proposition 3.2.} \textit{For the adjoint difference equation in Eq.
(\ref{42a}), it holds}
\begin{align}
&  \sigma(z)=\sigma^{\ast}(z-1)+\tau_{\nu-\mu}^{\ast}(z-1)\nabla x_{\nu-\mu
-1}(z),\label{44a}\\
&  \tau_{\nu-\mu}(z)=\frac{\sigma^{\ast}(z+1)-\sigma_{\nu-\mu}^{\ast
}(z-1)-\tau_{\nu-\mu}^{\ast}(z-1)\nabla x_{\nu-\mu-1}(z)}{\Delta x_{\nu-\mu
-1}(z)},\label{44b}\\
&  \lambda=\lambda_{\nu-\mu}^{\ast}-\Delta_{\nu-\mu-1}\left(  \frac{\tau
_{\nu-\mu}^{\ast}(z-1)\Delta x_{\nu-\mu-1}(z)-\nabla\sigma^{\ast}(z)}%
{\nabla_{\nu-\mu}x(z)}\right)  . \label{44c}%
\end{align}

\textbf{Proof. }From Eq. (\ref{41b}), we have%
\begin{equation}
\tau_{\nu-\mu}^{\ast}(z)\Delta x_{\nu-\mu-1}(z)=\sigma(z+1)-\sigma
(z-1)-\tau_{\nu-\mu}(z-1)\nabla x_{\nu-\mu-1}(z). \label{45}%
\end{equation}
Based on Eq. (\ref{41a}) and Eq. (\ref{45}), we obtain%
\[
\sigma(z+1)=\sigma^{\ast}(z)+\tau_{\nu-\mu}^{\ast}(z)\Delta x_{\nu-\mu-1}(z),
\]
implying that%
\[
\sigma(z)=\sigma^{\ast}(z-1)+\tau_{\nu-\mu}^{\ast}(z-1)\nabla x_{\nu-\mu
-1}(z).
\]
Thus, we obtain from Eq. (\ref{41a}) that%
\begin{align*}
\tau_{\nu-\mu}(z-1)  &  =\frac{{\sigma^{\ast}(z)-\sigma(z-1)}}{{\nabla
x_{\nu-\mu-1}(z)}}\\
&  =\frac{{\sigma^{\ast}(z)-\sigma^{\ast}(z-2)-\tau_{\nu-\mu}^{\ast
}(z-2)\nabla x_{\nu-\mu-1}(z-1)}}{{\nabla x_{\nu-\mu-1}(z)}},
\end{align*}
which is Eq. (\ref{44b}). Moreover, we obtain%
\begin{align}
\tau_{\nu-\mu}(z-1)\nabla x_{\nu-\mu-1}(z)-\nabla\sigma(z)=  &  \sigma^{\ast
}(z)-\sigma(z-1)-\nabla\sigma(z)\nonumber\\
=  &  \sigma^{\ast}(z)-\sigma(z)\nonumber\\
=  &  \sigma^{\ast}(z)-[\sigma^{\ast}(z-1)+\tau_{\nu-\mu}^{\ast}(z-1)\nabla
x_{\nu-\mu-1}(z)]\nonumber\\
=  &  \nabla\sigma^{\ast}(z)-\tau_{\nu-\mu}^{\ast}(z-1)\nabla x_{\nu-\mu
-1}(z). \label{46}%
\end{align}
Using Eq. (\ref{41c}) together with Eq. (\ref{46}) gives Eq. (\ref{44c}) and
hence the proof is completed.

Parallel to Corollary 3.1, one may obtain the following corollary.

\textbf{Corollary 3.2.} \textit{Eq. (\ref{44b}) and Eq. (\ref{44c}) can be
simplified to, respectively, }
\begin{equation}
\tau_{\nu-\mu}(z)=-\tau_{\nu-\mu-2}^{\ast}(z+1), \label{47a}%
\end{equation}%
\begin{equation}
\lambda=\lambda_{\nu-\mu}^{\ast}-\kappa_{2\nu-2\mu-1}^{\ast}. \label{47b}%
\end{equation}

\textbf{Proposition 3.3.} \textit{The adjoint difference equation Eq.
(\ref{42a}) can be rewritten as}
\begin{equation}
\sigma(z+1)\Delta_{\nu-\mu-1}\nabla_{\nu-\mu}w(z)-\tau_{\nu-\mu-2}%
(z+1)\nabla_{\nu-\mu}w(z)+(\lambda-\kappa_{2\nu-2\mu-1})w(z)=0. \label{48}%
\end{equation}

\textbf{Proof.} Since%
\[
\Delta_{\nu-\mu}w(z)-\nabla_{\nu-\mu}w(z)=\Delta(\frac{{\nabla w(z)}}%
{{\nabla_{\nu-\mu}x(z)}}),
\]
we have
\begin{align}
\tau_{\nu-\mu}^{\ast}(z)\Delta_{\nu-\mu}w(z)=  &  \tau_{\nu-\mu}^{\ast
}(z)\nabla_{\nu-\mu}w(z)+\tau_{\nu-\mu}^{\ast}(z)\Delta(\frac{{\nabla w(z)}%
}{{\nabla_{\nu-\mu}x(s)}})\nonumber\\
=  &  \tau_{\nu-\mu}^{\ast}(z)\nabla_{\nu-\mu}w(z)+\tau_{\nu-\mu}^{\ast
}(z)\Delta x_{\nu-\mu-1}(z)\frac{\Delta}{{\Delta x_{\nu-\mu-1}(z)}}%
(\frac{{\nabla w(z)}}{{\nabla_{\nu-\mu}x(z)}}). \label{49}%
\end{align}
Substituting Eq. (\ref{49}) into Eq. (\ref{42a}), we obtain%
\begin{align}
\hspace{-0.2in}  &  [\sigma^{\ast}(z)+\tau_{\nu-\mu}^{\ast}(z)\Delta
x_{\nu-\mu-1}(z)]\frac{\Delta}{{\Delta x_{\nu-\mu-1}(z)}}(\frac{{\nabla w(z)}%
}{{\nabla_{\nu-\mu}x(z)}})\nonumber\\
&  +\tau_{\nu-\mu}^{\ast}(z)\nabla_{\nu-\mu}w(z)+\lambda_{\nu-\mu}^{\ast
}w(z)\nonumber\\
=  &  0. \label{50}%
\end{align}
From Eq. (\ref{43a}), one may have
\begin{equation}
\sigma^{\ast}(z)+\tau_{\nu-\mu}^{\ast}(z)\Delta x_{\nu-\mu-1}(z)=\sigma(z+1).
\label{51}%
\end{equation}
Substituting it into Eq. (\ref{50}) and then using Eq. (\ref{43a}), we obtain
Eq. (\ref{48}) and hence complete the proof.

In the end of this section, we would like to prove that the adjoint difference
equation in Eq. (\ref{42a}) or Eq. (\ref{48}) is also a difference equation of
hypergeometric type on non-uniform lattices. To this end, we need only to
prove that
\[
\tilde{\sigma}^{\ast}(z)=\sigma^{\ast}(z)+\frac{1}{2}\tau_{\nu-\mu}^{\ast
}(z)\Delta x_{\nu-\mu-1}(z)=\sigma(z+1)+\frac{1}{2}\tau_{\nu-\mu-2}(z+1)\Delta
x_{\nu-\mu-1}(z)
\]
is a polynomial of degree at most two in the variable $x_{\nu-\mu}(z)$. In
fact, from Eq. (\ref{22}) and Lemma 2.1, we see that
\[
\tilde{\sigma}^{\ast}(z)=\sigma(z+1)+\frac{1}{2}\tau_{\nu-\mu-2}(z+1)\nabla
x_{\nu-\mu-1}(z+1)=\tilde{\sigma}_{\nu-\mu-2}(z+1)
\]
is a polynomial of degree at most two in the variable $x_{\nu-\mu
-2}(z+1)=x_{\nu-\mu}(z)$. Thus, we obtain the following theorem.

\textbf{Theorem 3.1.} \textit{The adjoint equation Eq. (\ref{48}) or}%
\begin{align}
&  \tilde{\sigma}_{\nu-\mu-2}(z+1)\Delta_{\nu-\mu-1}\nabla_{\nu-\mu
}w(z)\nonumber\\
&  -\frac{1}{2}\tau_{\nu-\mu-2}(z+1)[\Delta_{\nu-\mu}w(z)+\nabla_{\nu-\mu
}w(z)]+(\lambda-\kappa_{2\nu-2\mu-1})w(z)\nonumber\\
&  =0 \label{52a}%
\end{align}
\textit{is also a difference equation of hypergeometric type on non-uniform
lattices.}

By letting $\mu=\nu$ in the above equation, we immediately obtain the
following corollary.

\textbf{Corollary 3.3.} \textit{The adjoint equation of Eq. (\ref{NUSeq}) or}%
\begin{align}
&  \tilde{\sigma}_{-2}(z+1)\Delta_{-1}\nabla_{0}w(z)\nonumber\\
&  -\frac{1}{2}\tau_{-2}(z+1)[\Delta_{0}w(z)+\nabla_{0}w(z)]+(\lambda
-\kappa_{-1})w(z)\nonumber\\
=  &  0 \label{52b}%
\end{align}
\textit{is also a difference equation of hypergeometric type on non-uniform
lattices.}

\section{{\protect\Large Particular Solutions for Adjoint Difference
Equations}}

In this section, we first derive the forms of particular solutions for a
difference equation of hypergeometric type on non-uniform lattices (see
Proposition 4.1) and then use it to obtain the forms of particular solutions
for the adjoint difference equation given in Eq. (\ref{42a}) or an alternative
equation in Eq. (\ref{48}) (see Theorem 4.1). By letting $\mu=\nu$, one may
obtain the forms of particular solutions for the adjoint difference equation
given in Eq. (\ref{52b}) (see Theorem 4.2).

\textbf{Proposition 4.1. }\textit{On those classes of non-uniform lattices}
$x=x(z)$, \textit{the difference equation of hypergeometric type on
non-uniform lattices}
\begin{equation}
\sigma(z)\frac{\Delta}{{\Delta x_{\nu-\mu-1}(z)}}(\frac{{\nabla y(z)}}{{\nabla
x_{\nu-\mu}(z)}})-\tau_{\nu-\mu}(z)\frac{{\nabla y(z)}}{{\nabla x_{\nu-\mu
}(z)}}+\tilde{\lambda}y(z)=0\quad(\tilde{\lambda}\in R)
\end{equation}
\textit{has particular solutions in the form of}
\begin{equation}
y(z)=\sum\limits_{s=a}^{b-1}{\frac{{\rho_{\nu}(s)\nabla x_{\nu+1}(s)}%
}{{[x_{\nu}(s)-x_{\nu}(z-1)]^{(\mu-1)}}}},
\end{equation}
\textit{and also in the form of}
\begin{equation}
y(z)=\oint_{C}{\frac{{\rho_{\nu}(s)\nabla x_{\nu+1}(s)ds}}{{[x_{\nu}%
(s)-x_{\nu}(z-1)]^{(\mu-1)}}}},
\end{equation}
\textit{where }$a,b$ are complex with the same imaginary parts, $C$ \textit{is
a contour in the complex }$s$\textit{-plane, and} $x_{\nu}(s)=x(s+\frac{1}%
{2})$, \textit{if}

i) \textit{functions} $\rho(z)$ \textit{and} $\rho_{\gamma}(z)$
\textit{satisfy}
\begin{equation}
\frac{\Delta}{{\nabla x_{1}(z)}}[\sigma(z)\rho(z)]=\tau(z)\rho(z),\hspace
{0.1in}\frac{\Delta}{{\nabla x_{\nu+1}(z)}}[\sigma(z)\rho_{\nu}(z)]=\tau_{\nu
}(z)\rho_{\nu}(z);
\end{equation}

ii) $\mu,\nu$ \textit{satisfy}
\begin{equation}
\tilde{\lambda}+\kappa_{2\nu-(\mu-1)}\gamma(\mu-1)=0;
\end{equation}

iii) \textit{difference derivatives of the functions calculated by}
\begin{equation}
\phi_{\nu\mu}(z)=\sum\limits_{s=a}^{b-1}{\frac{{\rho_{\nu}(s)\nabla x_{\nu
-1}(s)}}{{[x_{\nu}(s)-x_{\nu}(z-1)]^{(\mu-1)}}}},
\end{equation}
\hspace{0.4in}\textit{or}
\begin{equation}
\phi_{\nu\mu}(z)=\oint_{C}{\frac{{\rho_{\nu}(s)\nabla x_{\nu-1}(s)ds}%
}{{[x_{\nu}(s)-x_{\nu}(z-1)]^{(\mu-1)}}}}%
\end{equation}
\hspace{0.4in}\textit{can be carried out by means of the formula}
\begin{equation}
\frac{{\nabla\phi_{\nu\mu}(z)}}{{\nabla x_{\nu-\mu}(z)}}=\gamma(\mu
-1)\phi_{\nu,\mu+1}(z);
\end{equation}

iv) \textit{the following equalities hold}\
\begin{equation}
\psi_{\nu\mu}(a,z)=\psi_{\nu\mu}(b,z),\hspace{0.1in}\oint_{C}{\Delta_{s}%
\psi_{\nu\mu}(s,z)}dz=0,
\end{equation}
\hspace{0.4in}\textit{where}
\begin{equation}
\psi_{\nu\mu}(s,z)=\frac{{\sigma(s)\rho_{\nu}(s)}}{{[x_{\nu-1}(s)-x_{\nu
-1}(z)]^{(\mu)}}}.
\end{equation}

\textbf{Proof.} To establish the relationship among $\Delta_{\nu-\mu-1}%
\nabla_{\nu-\mu}y(z),\nabla_{\nu-\mu}y(z),$ and $y(z)$, we need to find
nonzero functions $A_{i}(z),i=1,2,3$, such that
\begin{equation}
A_{1}(z)\Delta_{\nu-\mu-1}\nabla_{\nu-\mu}y(z)+A_{2}(z)\nabla_{\nu-\mu
}y(z)+A_{3}(z)y(z)=0. \label{63}%
\end{equation}
Note that
\[
\nabla_{\nu-\mu}y(z)=\gamma(\mu-1)\sum\limits_{s=a}^{b-1}{\frac{{\rho_{\nu
}(s)\nabla x_{\nu+1}(s)}}{{[x_{\nu}(s)-x_{\nu}(z-1)]^{(\mu)}}}},
\]%
\[
\Delta_{\nu-\mu-1}\nabla_{\nu-\mu}y(z)=\gamma(\mu)\gamma(\mu-1)\sum
\limits_{s=a}^{b-1}{\frac{{\rho_{\nu}(s)\nabla x_{\nu+1}(s)}}{{[x_{\nu
}(s)-x_{\nu}(z)]^{(\mu-1)}}}}.
\]
Substituting them into Eq. (\ref{63}) gives
\begin{align*}
\hspace{-0.2in}  &  \sum\limits_{s=a}^{b-1}{\frac{{\rho_{\nu}(s)\nabla
x_{\nu+1}(s)}}{{[x_{\nu}(s)-x_{\nu}(z)]^{(\mu+1)}}}}\{\gamma(\mu)\gamma
(\mu-1)A_{1}(z)+\gamma(\mu-1)A_{2}(z)[x_{\nu}(s)-x_{\nu}(z)]\\
&  +A_{3}(z)[x_{\nu}(s)-x_{\nu}(z)][x_{\nu}(s)-x_{\nu}(z-\mu)]\}\\
&  =\sum\limits_{s=a}^{b-1}{\frac{{\rho_{\nu}(s)\nabla x_{\nu+1}(s)}}%
{{[x_{\nu}(s)-x_{\nu}(z)]^{(\mu+1)}}}}P(s),
\end{align*}
where%
\begin{align*}
P(s)  &  =\gamma(\mu)\gamma(\mu-1)A_{1}(z)+\gamma(\mu-1)A_{2}(z)[x_{\nu
}(s)-x_{\nu}(z)]\\
&  +A_{3}(z)[x_{\nu}(s)-x_{\nu}(z)][x_{\nu}(s)-x_{\nu}(z-\mu)].
\end{align*}
On the other hand, we set%
\begin{align*}
\hspace{-0.2in}  &  \sum\limits_{s=a}^{b-1}{\frac{{\rho_{\nu}(s)\nabla
x_{\nu+1}(s)}}{{[x_{\nu}(s)-x_{\nu}(z)]^{(\mu+1)}}}}P(s)\\
&  =\sum\limits_{s=a}^{b-1}{\Delta_{s}[\frac{{\sigma(s)\rho_{\nu}(s)}%
}{{[x_{\nu-1}(s)-x_{\nu-1}(z)]^{(\mu)}}}}]\\
&  =\sum\limits_{s=a}^{b-1}{\frac{{\tau_{\nu}(s)\rho_{\nu}(s)\nabla x_{\nu
+1}(s)}}{{[x_{\nu-1}(s+1)-x_{\nu-1}(z)]^{(\mu)}}}}-\sum\limits_{s=a}%
^{b-1}{\frac{{\gamma(\mu)\sigma(s)\rho_{\nu}(s)\nabla x_{\nu+1}(s)}}{{[x_{\nu
}(s)-x_{\nu}(z)]^{(\mu+1)}}}}\\
&  =-\sum\limits_{s=a}^{b-1}{\frac{{\rho_{\nu}(s)\nabla x_{\nu+1}(s)}%
}{{[x_{\nu}(s)-x_{\nu}(z)]^{(\mu+1)}}}\{}\gamma(\mu)\sigma(s)-\tau_{\nu
}(s)[x_{\nu-\mu}(s)-x_{\nu-\mu}(z)]\}\\
&  =-\sum\limits_{s=a}^{b-1}{\frac{{\rho_{\nu}(s)\nabla x_{\nu+1}(s)}%
}{{[x_{\nu}(s)-x_{\nu}(z)]^{(\mu+1)}}}Q(s),}%
\end{align*}
where
\[
Q(s)=\gamma(\mu)\sigma(s)-\tau_{\nu}(s)[x_{\nu-\mu}(s)-x_{\nu-\mu}(z)].
\]
By Lemma 2.2, we have
\begin{align*}
Q(s)  &  =\gamma(\mu)\sigma(z)-\tau_{\nu-\mu}(z)[x_{\nu}(s)-x_{\nu}(z)]\\
&  -\kappa_{2\nu-(\mu-1)}[x_{\nu}(s)-x_{\nu}(z)][x_{\nu}(s)-x_{\nu}(z-\mu)].
\end{align*}
Comparing $P(s)$ with $Q(s)$ gives
\[
A_{1}(z)=\frac{1}{{\gamma(\mu-1)}}\sigma(z),\hspace{0.1in}A_{2}(z)=-\frac
{{\tau_{\nu-\mu}(z)}}{{\gamma(\mu-1)}},\hspace{0.1in}A_{3}(z)=-\kappa
_{2\nu-(\mu-1)},
\]
and hence the proof is completed .

\textbf{Theorem 4.1.} \textit{On those classes of non-uniform lattices}
$x=x(z)$, \textit{the adjoint difference equation given in Eq. (\ref{42a}) or
an alternative equation in Eq. (\ref{48}) as}
\begin{equation}
\sigma(z+1)\frac{\Delta}{{\Delta x_{\nu-\mu-1}(z)}}(\frac{{\nabla y(z)}%
}{{\nabla x_{\nu-\mu}(z)}})-\tau_{\nu-\mu-2}(z+1)\frac{{\nabla y(z)}}{{\nabla
x_{\nu-\mu}(z)}}+\lambda_{\nu-\mu}^{\ast}y(z)=0,
\end{equation}
\textit{has particular solutions in the form of}
\[
y(z)=\sum\limits_{s=a}^{b-1}{\frac{{\rho_{\nu}(s)\nabla x_{\nu+1}(s)}%
}{{[x_{\nu}(s)-x_{\nu}(z)]^{(\mu+1)}}}},
\]
\textit{and also in the form of}
\[
y(z)=\oint_{C}{\frac{{\rho_{\nu}(s)\nabla x_{\nu+1}(s)ds}}{{[x_{\nu}%
(s)-x_{\nu}(z)]^{(\mu+1)}}}},
\]
\textit{where} $C$ \textit{is a contour in the complex} $s$\textit{-plane,
and} $x_{\nu}(s)=x(s+\frac{1}{2})$, \textit{if}

i) \textit{functions} $\rho(z)$ \textit{and} $\rho_{\gamma}(z)$
\textit{satisfy}
\begin{equation}
\frac{\Delta}{{\nabla x_{1}(z)}}[\sigma(z)\rho(z)]=\tau(z)\rho(z),\hspace
{0.1in}\frac{\Delta}{{\nabla x_{\nu+1}(z)}}[\sigma(z)\rho_{\nu}(z)]=\tau_{\nu
}(z)\rho_{\nu}(z); \label{65}%
\end{equation}

ii) $\mu,\nu$ \textit{satisfy}
\begin{equation}
\lambda_{\nu-\mu}^{\ast}+\kappa_{2\nu-(\mu+1)}\gamma(\mu+1)=0;
\end{equation}

iii) \textit{difference derivatives of the functions calculated by}
\begin{equation}
\phi_{\nu\mu}(z)=\sum\limits_{s=a}^{b-1}{\frac{{\rho_{\nu}(s)\nabla x_{\nu
-1}(s)}}{{[x_{\nu}(s)-x_{\nu}(z)]^{(\mu+1)}}}},
\end{equation}
\hspace{0.4in}\textit{or}
\begin{equation}
\phi_{\nu\mu}(z)=\oint_{C}{\frac{{\rho_{\nu}(s)\nabla x_{\nu-1}(s)ds}%
}{{[x_{\nu}(s)-x_{\nu}(z)]^{(\mu+1)}}}}%
\end{equation}
\hspace{0.4in}\textit{can be carried out by means of the formula}
\begin{equation}
\frac{{\nabla\phi_{\nu\mu}(z)}}{{\nabla x_{\nu-\mu}(z)}}=\gamma(\mu
+1)\phi_{\nu,\mu+1}(z);
\end{equation}

iv) \textit{the following equalities hold}\
\begin{equation}
\psi_{\nu\mu}(a,z)=\psi_{\nu\mu}(b,z),\hspace{0.1in}\oint_{C}{\Delta_{s}%
\psi_{\nu\mu}(s,z)}dz=0
\end{equation}
\hspace{0.4in}\textit{where}
\begin{equation}
\psi_{\nu\mu}(s,z)=\frac{{\sigma(s)\rho_{\nu}(s)}}{{[x_{\nu-1}(s)-x_{\nu
-1}(z+1)]^{(\mu+1)}}}.
\end{equation}

\textbf{Proof.} Eq. (\ref{65}) can be written as%
\begin{align*}
&  \hspace{-0.2in}\sigma(z+1)\frac{\Delta}{{\Delta x_{\nu-\mu-3}(z+1)}}%
(\frac{{\nabla y(z)}}{{\nabla x_{\nu-\mu-2}(z+1)}})-\tau_{\nu-\mu-2}%
(z+1)\frac{{\nabla y(z)}}{{\nabla x_{\nu-\mu-1}(z+1)}}\\
&  -\kappa_{2\nu-(\mu+1)}\gamma(\mu+1)y(z)\\
&  =0.
\end{align*}
Letting $y(z)=\tilde{y}(z+1)$, we obtain
\begin{align*}
&  \hspace{-0.2in}\sigma(z)\frac{\Delta}{{\Delta x_{\nu-\mu-3}(z)}}%
(\frac{{\nabla\tilde{y}(z)}}{{\nabla x_{\nu-\mu-2}(z)}})-\tau_{\nu-\mu
-2}(z)\frac{{\nabla\tilde{y}(z)}}{{\nabla x_{\nu-\mu-2}(z)}}-\kappa_{2\nu
-(\mu+1)}\gamma(\mu+1)\tilde{y}(z)\\
&  =0.
\end{align*}
Letting $\mu+2=\tilde{\mu}$, we further obtain
\[
\sigma(z)\frac{\Delta}{{\Delta x_{\nu-\tilde{\mu}-1}(z)}}(\frac{{\nabla
\tilde{y}(z)}}{{\nabla x_{\nu-\tilde{\mu}}(z)}})-\tau_{\nu-\tilde{\mu}%
}(z)\frac{{\nabla\tilde{y}(z)}}{{\nabla x_{\nu-\tilde{\mu}}(z)}}-\kappa
_{2\nu-(\tilde{\mu}-1)}\gamma(\tilde{\mu}-1)\tilde{y}(z)=0.
\]
By Proposition 4.1, we have
\[
y(z)=\tilde{y}(z+1)=\sum\limits_{s=a}^{b-1}{\frac{{\rho_{\nu}(s)\nabla
x_{\nu+1}(s)}}{{[x_{\nu}(s)-x_{\nu}(z)]^{(\tilde{\mu}-1)}}}=\sum
\limits_{s=a}^{b-1}{\frac{{\rho_{\nu}(s)\nabla x_{\nu+1}(s)}}{{[x_{\nu
}(s)-x_{\nu}(z)]^{(\mu+1)}}}}}%
\]
and
\[
y(z)=\tilde{y}(z+1)=\oint_{C}{\frac{{\rho_{\nu}(s)\nabla x_{\nu+1}(s)ds}%
}{{[x_{\nu}(s)-x_{\nu}(z)]^{(\tilde{\mu}-1)}}}}=\oint_{C}{\frac{{\rho_{\nu
}(s)\nabla x_{\nu+1}(s)ds}}{{[x_{\nu}(s)-x_{\nu}(z)]^{(\mu+1)}}}},
\]
and hence complete the proof.

Finally, by letting $\mu=\nu$, then from Eq. (\ref{43b}), $\lambda^{\ast
}=\lambda-\kappa_{-1},$ one may obtain the forms of particular solutions for
the adjoint difference equation in \textit{Eq. (\ref{52b}) }.

\textbf{Theorem 4.2.} \textit{On those classes of non-uniform lattices}
$x=x(z)$, \textit{the adjoint difference equation given in Eq. (\ref{52b})
as}
\begin{equation}
\sigma(z+1)\frac{\Delta}{{\Delta x_{-1}(z)}}(\frac{{\nabla y(z)}}{{\nabla
x(z)}})-\tau_{-2}(z+1)\frac{{\nabla y(z)}}{{\nabla x(z)}}+\lambda^{\ast
}y(z)=0,
\end{equation}
\textit{has particular solutions in the form of}
\[
y(z)=\sum\limits_{s=a}^{b-1}{\frac{{\rho_{\nu}(s)\nabla x_{\nu+1}(s)}%
}{{[x_{\nu}(s)-x_{\nu}(z)]^{(\nu+1)}}}},
\]
\textit{and also in the form of}
\[
y(z)=\oint_{C}{\frac{{\rho_{\nu}(s)\nabla x_{\nu+1}(s)ds}}{{[x_{\nu}%
(s)-x_{\nu}(z)]^{(\nu+1)}}}},
\]
\textit{where} $C$ \textit{is a contour in the complex} $s$\textit{-plane,
and} $x_{\nu}(s)=x(s+\frac{1}{2})$, \textit{if}

i) \textit{functions} $\rho(z)$ \textit{and} $\rho_{\gamma}(z)$
\textit{satisfy}
\begin{equation}
\frac{\Delta}{{\nabla x_{1}(z)}}[\sigma(z)\rho(z)]=\tau(z)\rho(z),\hspace
{0.1in}\frac{\Delta}{{\nabla x_{\nu+1}(z)}}[\sigma(z)\rho_{\nu}(z)]=\tau_{\nu
}(z)\rho_{\nu}(z); \label{66}%
\end{equation}

ii) $\mu,\nu$ \textit{satisfy}
\begin{equation}
\lambda^{\ast}+\kappa_{\nu-1}\gamma(\nu+1)=0;
\end{equation}

iii) \textit{difference derivatives of the functions calculated by}
\begin{equation}
\phi_{\nu\nu}(z)=\sum\limits_{s=a}^{b-1}{\frac{{\rho_{\nu}(s)\nabla x_{\nu
-1}(s)}}{{[x_{\nu}(s)-x_{\nu}(z)]^{(\nu+1)}}}},
\end{equation}
\hspace{0.4in}\textit{or}
\begin{equation}
\phi_{\nu\nu}(z)=\oint_{C}{\frac{{\rho_{\nu}(s)\nabla x_{\nu-1}(s)ds}%
}{{[x_{\nu}(s)-x_{\nu}(z)]^{(\nu+1)}}}}%
\end{equation}
\hspace{0.4in}\textit{can be carried out by means of the formula}
\begin{equation}
\frac{{\nabla\phi_{\nu\nu}(z)}}{{\nabla x(z)}}=\gamma(\nu+1)\phi_{\nu,\nu
+1}(z);
\end{equation}

iv) \textit{the following equalities hold}\
\begin{equation}
\psi_{\nu\nu}(a,z)=\psi_{\nu\nu}(b,z),\hspace{0.1in}\oint_{C}{\Delta_{s}%
\psi_{\nu\nu}(s,z)}dz=0
\end{equation}
\hspace{0.4in}\textit{where}
\begin{equation}
\psi_{\nu\nu}(s,z)=\frac{{\sigma(s)\rho_{\nu}(s)}}{{[x_{\nu-1}(s)-x_{\nu
-1}(z+1)]^{(\nu+1)}}}.
\end{equation}

\section{{\protect\Large Application and New Fundamental Theorems}}

Based on Proposition 3.1 and Theorem 4.1, one may obtain the following corollary.

\textbf{Corollary 5.1.} \textit{Under the hypotheses of Theorem 4.1, the
equation}%
\begin{equation}
\sigma(z)\frac{\Delta}{{\Delta x_{\nu-\mu-1}(z)}}(\frac{{\nabla y(z)}}{{\nabla
x_{\nu-\mu}(z)}})+\tau_{\nu-\mu}(z)\frac{{\Delta y(z)}}{{\Delta x_{\nu-\mu
}(z)}}+\lambda y(z)=0 \label{72}%
\end{equation}
\textit{has particular solutions in the form of}%
\begin{equation}
y(z)=\frac{1}{{\rho_{\nu-\mu}(z)}}\sum\limits_{s=a}^{b-1}{\frac{{\rho_{\nu
}(s)\nabla x_{\nu+1}(s)}}{{[x_{\nu}(s)-x_{\nu}(z)]^{(\mu+1)}}}}%
\end{equation}
\textit{and also the form of}%
\begin{equation}
y(z)=\frac{1}{{\rho_{\nu-\mu}(z)}}\oint_{C}{\frac{{\rho_{\nu}(s)\nabla
x_{\nu+1}(s)ds}}{{[x_{\nu}(s)-x_{\nu}(z)]^{(\mu+1)}}}},
\end{equation}
\textit{where} $\rho(z),\rho_{\nu}(z)$ \textit{satisfy}%
\[
\frac{{\Delta(\sigma(z)\rho(z))}}{{\nabla x_{1}(z)}}=\tau(z)\rho
(z),\hspace{0.1in}\frac{{\Delta(\sigma(z)\rho_{\nu}(z))}}{{\nabla x_{\nu
+1}(z)}}=\tau_{\nu}(z)\rho_{\nu}(z),
\]
\textit{and} $\nu,\mu$ \textit{are roots of the equation}%
\[
\lambda+\kappa_{2\nu-\mu}\gamma(\mu)=0.
\]

Note that by letting $\mu=\nu$ in Corollary 5.1, Eq. (\ref{72}) can be reduced
to Eq. (\ref{NUSeq}). Thus, we obtain the following well-known theorem given
in \cite{suslov1989}.

\textbf{Corollary 5.2 }(Theorem 2.2 in \cite{suslov1989})\textbf{.}
\textit{Under the hypotheses of Corollary 5.1 with} $\mu=\nu$, \textit{the
equation}%
\begin{equation}
\sigma(z)\frac{\Delta}{{\Delta x_{-1}(z)}}(\frac{{\nabla y(z)}}{{\nabla
x_{0}(z)}})+\tau(z)\frac{{\Delta y(z)}}{{\Delta x_{0}(z)}}+\lambda y(z)=0
\end{equation}
\textit{has particular solutions in the form of}
\begin{equation}
y(z)=\frac{1}{{\rho(z)}}\sum\limits_{s=a}^{b-1}{\frac{{\rho_{\nu}(s)\nabla
x_{\nu+1}(s)}}{{[x_{\nu}(s)-x_{\nu}(z)]^{(\nu+1)}}}}%
\end{equation}
\textit{and also in the form of}
\begin{equation}
y(z)=\frac{1}{{\rho(z)}}\oint_{C}{\frac{{\rho_{\nu}(s)\nabla x_{\nu+1}(s)ds}%
}{{[x_{\nu}(s)-x_{\nu}(z)]^{(\nu+1)}}}},
\end{equation}
\textit{where} $\rho(z),\rho_{\nu}(z)$ \textit{satisfy}
\begin{equation}
\frac{{\Delta(\sigma(z)\rho(z))}}{{\nabla x_{1}(z)}}=\tau(z)\rho
(z),\hspace{0.1in}\frac{{\Delta(\sigma(z)\rho_{\nu}(z))}}{{\nabla x_{\nu
+1}(z)}}=\tau_{\nu}(z)\rho_{\nu}(z), \label{pearson1}%
\end{equation}
\textit{and} $\nu$ \textit{is the root of the equation}
\[
\lambda+\kappa_{\nu}\gamma(\nu)=0.
\]

\textbf{Remark 5.1.} It should be pointed out that Theorem 4.1 may be obtained
based on the Suslov theorem coupled with Proposition 3.1. However, without
using Proposition 3.1, Theorem 4.1 seems cannot be obtained simply using the
Suslov theorem without coupling with Proposition 3.1. We consider Theorem 4.1
to be a new result because we prove it directly and have not seen it as well
as Proposition 3.1 from the other literatures. Reversely, the Suslov theorem
(Corollary 5.2) can be obtained based on Theorem 4.1 and Proposition 3.1. This
new proof not only gives another way to prove the Suslov theorem, but also is
our purpose showing the important application of the obtained adjoint
equations and their solutions in this study.

\textbf{Remark 5.2. }One of\textbf{ }interests for the adjoint equation in
\cite{area2003,area2005} is because it can be used to find the general
solution of the hypergeometric and $q-$hypergeometric equation. Using our
results obtained in this study, it is possible to do something similar.
Indeed, we have done some work where the results can be seen in our recent
manuscript \cite{cheng2018}.

We now prove another kind of fundamental theorems for Eq. (\ref{NUSeq}) and
Eq. (\ref{72}), respectively, which are essentially new results and their
expressions are different from the Suslov theorem (as seen in Corollaries 5.1
and 5.2).

\textbf{Theorem 5.1. }\textit{On those classes of nonuniform lattices
}$x=x(z)$, \textit{the difference equation of hypergeometric type on
non-uniform lattices}%
\begin{equation}
\sigma(z)\frac{\Delta}{{\Delta x_{\nu-\mu-1}(z)}}(\frac{{\nabla y(z)}}{{\nabla
x_{\nu-\mu}(z)}})+\tau_{\nu-\mu}(z)\frac{{\Delta y(z)}}{{\Delta x_{\nu-\mu
}(z)}}+\lambda y(z)=0 \label{78}%
\end{equation}
\textit{has particular solutions in the form of}
\begin{equation}
y(z)=\sum\limits_{s=a}^{b-1}{{{[x_{\nu}(s)-x_{\nu}(z)]^{(\mu+1)}}}}{{\rho
_{\nu}(s)\nabla x_{\nu+1}(s)}},
\end{equation}
\textit{and also in the form of}%
\begin{equation}
y(z)=\oint_{C}{{{[x_{\nu}(s)-x_{\nu}(z)]^{(\mu+1)}}}}{{\rho_{\nu}(s)\nabla
x_{\nu+1}(s)ds}},
\end{equation}
\textit{where} $C$ \textit{is a contour in the complex }$s$\textit{-plane}%
,\textit{\ and} $x_{\nu}(s)=x(s+\frac{1}{2})$, \textit{if}

i) \textit{functions} $\rho(z)$ \textit{and} $\rho_{\gamma}(z)$
\textit{satisfy}%
\begin{equation}
\frac{\nabla}{{\nabla x_{\nu+1}(z)}}[\sigma(z)\rho_{\nu}(z)]+\tau_{\nu}%
(z)\rho_{\nu}(z)=0;
\end{equation}

ii) $\mu,\nu$ \textit{satisfy the equation}
\begin{equation}
\lambda+\kappa_{2\nu-(\mu-1)}\gamma(\mu+1)=0;
\end{equation}

ii) \textit{difference derivatives of the functions calculated by}
\begin{equation}
\phi_{\nu\mu}(z)=\sum\limits_{s=a}^{b-1}{{{[x_{\nu}(s)-x_{\nu}(z)]^{(\mu+1)}}%
}}{{\rho_{\nu}(s)\nabla x_{\nu-1}(s)}},
\end{equation}
\hspace{0.4in}\textit{or}
\begin{equation}
\phi_{\nu\mu}(z)=\oint_{C}{{{[x_{\nu}(s)-x_{\nu}(z)]^{(\mu+1)}}}}{{\rho_{\nu
}(s)\nabla x_{\nu-1}(s)ds}}%
\end{equation}
\hspace{0.4in}\textit{can be carried out by means of the formula}
\begin{equation}
\frac{{\Delta\phi_{\nu\mu}(z)}}{{\Delta x_{\nu-\mu}(z)}}=-\gamma(\mu
+1)\phi_{\nu,\mu-1}(z);
\end{equation}

iv) \textit{the following equalities hold}\
\begin{equation}
\psi_{\nu\mu}(a,z)=\psi_{\nu\mu}(b,z),\hspace{0.1in}\oint_{C}{\nabla_{s}%
\psi_{\nu\mu}(s,z)}dz=0,
\end{equation}
\hspace{0.4in}\textit{where}
\begin{equation}
\psi_{\nu\mu}(s,z)={\sigma(s)\rho_{\nu}(s)}{{[x_{\nu+1}(s)-x_{\nu
+1}(z-1)]^{(\mu)}}}.
\end{equation}

\textbf{Proof.} To establish the relationship among $\Delta_{\nu-\mu-1}%
\nabla_{\nu-\mu}y(z),\Delta_{\nu-\mu}y(z),$ and $y(z)$, we need to find
nonzero functions $A_{i}(z),i=1,2,3$, such that
\[
A_{1}(z)\Delta_{\nu-\mu-1}\nabla_{\nu-\mu}y(z)+A_{2}(z)\nabla_{\nu-\mu
}y(z)+A_{3}(z)y(z)=0.
\]
Substituting%
\[
\Delta_{\nu-\mu}y(z)=-\gamma(\mu+1)\sum\limits_{s=a}^{b-1}{{{[x_{\nu
}(s)-x_{\nu}(z)]^{(\mu)}}}}{{\rho_{\nu}(s)\nabla x_{\nu+1}(s)}},
\]%
\[
\Delta_{\nu-\mu-1}\nabla_{\nu-\mu}y(z)=\gamma(\mu+1)\gamma(\mu)\sum
\limits_{s=a}^{b-1}{{{[x_{\nu}(s)-x_{\nu}(z-1)]^{(\mu-1)}}}}{{\rho_{\nu
}(s)\nabla x_{\nu+1}(s)}}%
\]
into Eq. (\ref{78}), we obtain%
\begin{align*}
\hspace{-0.2in}  &  \sum\limits_{s=a}^{b-1}{{{\rho_{\nu}(s)\nabla x_{\nu
+1}(s)}}{{[x_{\nu}(s)-x_{\nu}(z-1)]^{(\mu-1)}}}}\{\gamma(\mu+1)\gamma
(\mu)A_{1}(z)\\
&  -\gamma(\mu+1)A_{2}(z)[x_{\nu}(s)-x_{\nu}(z)]+A_{3}(z)[x_{\nu}(s)-x_{\nu
}(z)][x_{\nu}(s)-x_{\nu}(z-\mu)]\}\\
&  =\sum\limits_{s=a}^{b-1}{{{\rho_{\nu}(s)\nabla x_{\nu+1}(s)}}{{[x_{\nu
}(s)-x_{\nu}(z-1)]^{(\mu-1)}}}}P(s),
\end{align*}
where
\begin{align*}
P(s)  &  =\{\gamma(\mu+1)\gamma(\mu)A_{1}(z)-\gamma(\mu+1)A_{2}(z)[x_{\nu
}(s)-x_{\nu}(z)]\\
&  +A_{3}(z)[x_{\nu}(s)-x_{\nu}(z)][x_{\nu}(s)-x_{\nu}(z-\mu)]\}.
\end{align*}
On the other hand, we let%
\begin{align*}
&  \sum\limits_{s=a}^{b-1}{{{\rho_{\nu}(s)\nabla x_{\nu+1}(s)}}{{[x_{\nu
}(s)-x_{\nu}(z-1)]^{(\mu-1)}}}}P(s)\\
&  =\sum\limits_{s=a}^{b-1}{\nabla_{s}[{\sigma(s)\rho_{\nu}(s)}{{[x_{\nu
+1}(s)-x_{\nu+1}(z-1)]^{(\mu)}}}}]\\
&  =-\sum\limits_{s=a}^{b-1}{{{\tau_{\nu}(s)\rho_{\nu}(s)\nabla x_{\nu+1}(s)}%
}{{[x_{\nu+1}(s-1)-x_{\nu+1}(z-1)]^{(\mu)}}}}\\
&  +\sum\limits_{s=a}^{b-1}{{{\gamma(\mu)\sigma(s)\rho_{\nu}(s)\nabla
x_{\nu+1}(s)}}{{[x_{\nu}(s)-x_{\nu}(z-1)]^{(\mu-1)}}}}\\
&  =\sum\limits_{s=a}^{b-1}{{{\rho_{\nu}(s)\nabla x_{\nu+1}(s)}}{{[x_{\nu
}(s)-x_{\nu}(z-1)]^{(\mu-1)}}}\{}\gamma(\mu)\sigma(s)-\tau_{\nu}(s)[x_{\nu
-\mu}(s)-x_{\nu-\mu}(z)]\}\\
&  =\sum\limits_{s=a}^{b-1}{{{\rho_{\nu}(s)\nabla x_{\nu+1}(s)}}{{[x_{\nu
}(s)-x_{\nu}(z-1)]^{(\mu-1)}}}Q(s),}%
\end{align*}
where
\[
Q(s)=\gamma(\mu)\sigma(s)-\tau_{\nu}(s)[x_{\nu-\mu}(s)-x_{\nu-\mu}(z)].
\]
By Lemma 2.2, we have
\begin{align*}
Q(s)  &  =\gamma(\mu)\sigma(z)-\tau_{\nu-\mu}(z)[x_{\nu}(s)-x_{\nu}(z)]\\
&  -\kappa_{2\nu-(\mu-1)}[x_{\nu}(s)-x_{\nu}(z)][x_{\nu}(s)-x_{\nu}(z-\mu)]
\end{align*}
Comparing $P(s)$ with $Q(s)$ gives
\[
A_{1}(z)=\frac{1}{{\gamma(\mu+1)}}\sigma(z),\hspace{0.1in}A_{2}(z)=\frac
{{\tau_{\nu-\mu}(z)}}{{\gamma(\mu+1)}},\hspace{0.1in}A_{3}(z)=-\kappa
_{2\nu-(\mu-1)},
\]
and hence we have completed the proof.

Letting $\mu=\nu$ in Theorem 5.1 gives the following theorem.

\textbf{Theorem 5.2. }\textit{Under the hypotheses of Theorem 5.1 with}
$\mu=\nu$, \textit{the equation}%
\begin{equation}
\sigma(z)\frac{\Delta}{{\Delta x_{-1}(z)}}(\frac{{\nabla y(z)}}{{\nabla
x_{0}(z)}})+\tau(z)\frac{{\Delta y(z)}}{{\Delta x_{0}(z)}}+\lambda y(z)=0
\label{th52}%
\end{equation}
\textit{has particular solutions in the form of}%
\begin{equation}
y(z)=\sum\limits_{s=a}^{b-1}{{[x_{\nu}(s)-x_{\nu}(z)]^{(\nu+1)}}}{{\rho_{\nu
}(s)\nabla x_{\nu+1}(s)}}%
\end{equation}
\textit{and also in the form of}
\begin{equation}
y(z)=\oint_{C}{{{[x_{\nu}(s)-x_{\nu}(z)]^{(\nu+1)}}}}{{\rho_{\nu}(s)\nabla
x_{\nu+1}(s)ds}}, \label{newformula}%
\end{equation}
\textit{where} $\rho(z),\rho_{\nu}(z)$ \textit{satisfy}
\begin{equation}
\frac{{\nabla(\sigma(z)\rho_{\nu}(z))}}{{\nabla x_{\nu+1}(z)}}+\tau_{\nu
}(z)\rho_{\nu}(z)=0, \label{pearson2}%
\end{equation}
\textit{and} $\nu$ \textit{is the root of the equation}
\begin{equation}
\lambda+\kappa_{\nu+1}\gamma(\nu+1)=0.
\end{equation}

In contrast with obtaining the solution $\rho_{\nu}(z)$ from the well-known
Pearson equation in \ Eq. (\ref{pearson1}), it seems more difficult to obtain
$\rho_{\nu}(z)$ directly from Eq. (\ref{pearson2}). However, coupling Eq.
(\ref{pearson2}) with Eq. (\ref{pearson1}), we may build up a useful
relationship between them as described in the following lemma.

\textbf{Lemma 5.1. }Let $\widetilde{{\rho_{\nu}}}{(z)},$ ${\rho_{\nu}(z)}$
satisfy the Pearson equation
\begin{equation}
\frac{{\Delta(\sigma(z)\widetilde{{\rho_{\nu}}}{(z)})}}{{\nabla x_{\nu+1}(z)}%
}=\tau_{\nu}(z)\widetilde{{\rho_{\nu}}}{(z)}, \label{pson1}%
\end{equation}
and%
\begin{equation}
\frac{{\nabla(\sigma(z)\rho_{\nu}(z))}}{{\nabla x_{\nu+1}(z)}}+\tau_{\nu
}(z)\rho_{\nu}(z)=0, \label{pson2}%
\end{equation}
then it holds%
\begin{equation}
{\rho_{\nu}(z)=}\frac{{const}}{{\sigma(z)\sigma(z+1)\widetilde{{\rho_{\nu}}%
}{(z+1)}}}. \label{lem}%
\end{equation}

\textbf{Proof. }Note that%
\begin{align*}
&  {\Delta\lbrack\sigma(z)\widetilde{{\rho_{\nu}}}{(z)}\sigma(z-1)\rho_{\nu
}(z-1)]}\\
&  ={\sigma(z)\widetilde{{\rho_{\nu}}}{(z)}\Delta\lbrack\sigma(z-1)\rho_{\nu
}(z-1)]+\sigma(z)\rho_{\nu}(z)\Delta\lbrack\sigma(z)\widetilde{{\rho_{\nu}}%
}{(z)}}\\
&  ={\sigma(z)\widetilde{{\rho_{\nu}}}{(z)}\nabla\lbrack\sigma(z)\rho_{\nu
}(z)]+\sigma(z)\rho_{\nu}(z)\Delta\lbrack\sigma(z)\widetilde{{\rho_{\nu}}%
}{(z).}}%
\end{align*}
Based on Eq. (\ref{pson1}) and Eq. (\ref{pson2}), one may obtain%
\begin{align*}
&  {\Delta\lbrack\sigma(z)\widetilde{{\rho_{\nu}}}{(z)}\sigma(z-1)\rho_{\nu
}(z-1)]}\\
&  =-{\sigma(z)\widetilde{{\rho_{\nu}}}{(z)}}\tau_{\nu}(z)\rho_{\nu}(z){\nabla
x_{\nu+1}(z)}+{\sigma(z)\rho_{\nu}(z)}\tau_{\nu}(z)\widetilde{{\rho_{\nu}}%
}{(z)\nabla x_{\nu+1}(z)}\\
&  =0,
\end{align*}
which yields%
\[
{\sigma(z)\widetilde{{\rho_{\nu}}}{(z)}\sigma(z-1)\rho_{\nu}(z-1)=const.}%
\]
Hence, Eq. (\ref{lem}) is obtained.

In particular, when the quadratic lattice\textbf{ }$x(z)=z^{2}$, we have the
following lemma.

\textbf{Lemma 5.2. }For $x(z)=z^{2}$, let ${\rho_{\nu}(z)}$ satisfy Eq.
(\ref{pson2}), then%
\begin{equation}
\frac{{\rho_{\nu}(z+1)}}{{\rho_{\nu}(z)}}=\frac{{\sigma(z)}}{{\sigma
(-z-1-\nu)}}. \label{lem2}%
\end{equation}

\textbf{Proof. }For $x(z)=z^{2},$ we have $\nabla x_{1}(z)=2z$ and the
property (also seen in Eq. (3.10.4) on page 123 in \cite{nikiforov1991})%
\begin{equation}
\sigma(z)+\tau(z)\nabla x_{1}(z)=\sigma(-z). \label{52}%
\end{equation}
Let ${\widetilde{{\rho_{\nu}}}{(z)}}$ satisfy Eq. (\ref{pson1}). From Eq.
(\ref{23}) and Eq. (\ref{52}), we obtain%
\begin{align}
\frac{{\widetilde{{\rho_{\nu}}}{(z+1)}}}{{\widetilde{{\rho_{\nu}}}{(z)}}}  &
=\frac{\sigma(z)+\tau_{\nu}(z)\nabla x_{\nu+1}(z)}{\sigma(z+1)}\nonumber\\
&  =\frac{\sigma(z+\nu)+\tau(z+\nu)\nabla x_{1}(z+\nu)}{\sigma(z+1)}%
\nonumber\\
&  =\frac{\sigma(-z-\nu)}{\sigma(z+1)}. \label{coro}%
\end{align}
By Lemma 5.1\textbf{ }and Eq. (\ref{coro}), we have%
\begin{align*}
\frac{{\rho_{\nu}(z+1)}}{{\rho_{\nu}(z)}}  &  =\frac{{\sigma(z)\sigma
(z+1)\widetilde{{\rho_{\nu}}}{(z+1)}}}{{\sigma(z+1)\sigma(z+2)\widetilde{{\rho
_{\nu}}}{(z+2)}}}\\
&  =\frac{{\sigma(z)}}{{\sigma(z+2)}}\frac{{\widetilde{{\rho_{\nu}}}{(z+1)}}%
}{{\widetilde{{\rho_{\nu}}}{(z+2)}}}\\
&  =\frac{{\sigma(z)}}{{\sigma(z+2)}}\frac{\sigma(z+2)}{\sigma(-z-1-\nu)}\\
&  =\frac{\sigma(z)}{\sigma(-z-1-\nu)},
\end{align*}
and complete the proof.

Finally, we give an example to illustrate the application of Theorem 5.2 for
the case of the quadratic lattice\textbf{ }$x(z)=z^{2}$.

\textbf{Example 5.1. }Consider the equation%
\[
\sigma(z)\frac{\Delta}{{\Delta x_{-1}(z)}}(\frac{{\nabla y(z)}}{{\nabla
x_{0}(z)}})+\tau(z)\frac{{\Delta y(z)}}{{\Delta x_{0}(z)}}+\lambda y(z)=0.
\]
where the lattice $x(s)=s^{2}$, $\sigma(s)=\prod\limits_{k=1}^{4}(s-s_{k}),$
and $s_{k},k=1,2,3,4$, are arbitrary complex numbers, which give
$\sigma(-s-1-\nu)=\prod\limits_{k=1}^{4}(s_{k}-s-\nu-1)$. We would like to
find its solution.

\textbf{Solution.} From Eq. (\ref{lem2}), we have%
\begin{equation}
\frac{{\rho_{\nu}(s+1)}}{{\rho_{\nu}(s)}}=\frac{\sigma(s)}{\sigma(-s-1-\nu
)}=\prod\limits_{k=1}^{4}\frac{(s-s_{k})}{(-s_{k}-s-\nu-1)}. \label{pson3}%
\end{equation}
Since
\[
\frac{(s-s_{k})}{(-s_{k}-s-\nu-1)}=\frac{\Gamma(s+1-s_{k})\Gamma(-s_{k}%
-s-\nu-1)}{\Gamma(s-s_{k})\Gamma(-s_{k}-s-\nu)},
\]
we choose a solution of Eq. (\ref{pson3}) in the form%
\begin{align*}
{\rho_{\nu}{(s)}}  &  {=}{C}_{0}\prod\limits_{k=1}^{4}\Gamma(s-s_{k}%
)\Gamma(-s_{k}-s-\nu)\sin2\pi(s+\frac{\nu+1}{2}),\\
C_{0}^{-1}  &  =\frac{\sin\pi(s-z+\nu+1)}{\sin\pi(s-z)}.
\end{align*}
Using the "generalized power" in the form given in \cite{suslov1989}%
\begin{align*}
{{{[x_{\nu}(s)-x_{\nu}(z)]^{(\nu+1)}}}}  &  {=}{{}}\frac{\Gamma(s-z+\nu
+1)\Gamma(s+z+\nu+1)}{\Gamma(s-z)\Gamma(s+z)},\\
x(z)  &  =z^{2},
\end{align*}
and%
\begin{align*}
\frac{\Gamma(s-z+\nu+1)}{\Gamma(s-z)}  &  =\frac{\pi/[\sin\pi(s-z+\nu
+1)\Gamma(z-s-\nu)]}{\pi/[\sin\pi(s-z)\Gamma(1+z-s)]}\\
&  =\frac{\sin\pi(s-z)\Gamma(1+z-s)}{\sin\pi(s-z+\nu+1)\Gamma(z-s-\nu)}\\
&  =C_{0}^{-1}\frac{\Gamma(1+z-s)}{\Gamma(z-s-\nu)},
\end{align*}
we obtain%
\begin{align*}
{{{[x_{\nu}(s)-x_{\nu}(z)]^{(\nu+1)}}}}  &  {=}{{}}C_{0}^{-1}\frac
{\Gamma(1+z-s)\Gamma(s+z+\nu+1)}{\Gamma(z-s-\nu)\Gamma(s+z)},\\
x(z)  &  =z^{2}.
\end{align*}
Based on Eq. (\ref{newformula}) in Theorem 5.2, we obtain%
\begin{align*}
y_{\nu}(z)  &  =\oint_{C}{{{[x_{\nu}(s)-x_{\nu}(z)]^{(\nu+1)}}}}{{\rho_{\nu
}(s)\nabla x_{\nu+1}(s)ds}}\\
&  =\oint_{C}\frac{\Gamma(1+z-s)\Gamma(s+z+\nu+1)}{\Gamma(z-s-\nu)\Gamma
(s+z)}{{\prod\limits_{k=1}^{4}\Gamma(s-s_{k})\Gamma(-s_{k}-s-\nu)(2s+\nu)ds.}}%
\end{align*}
Setting ${{2s+\nu=2t,}}$ we obtain%
\begin{align*}
2t  &  =\frac{\Gamma(1+2t)}{\Gamma(2t)}\\
&  =\frac{\pi}{\Gamma(2t)\Gamma(-2t)\sin\pi(-2t)}\\
&  =\frac{\pi}{\Gamma(2t)\Gamma(-2t)\sin2\pi(t+\frac{1}{2})}.
\end{align*}
Thus, we obtain a solution%
\begin{align*}
y_{\nu}(z)  &  =\pi\int_{-i\infty}^{i\infty}\frac{\Gamma(1+z+\frac{\nu}%
{2}-t)\Gamma(1+z+\frac{\nu}{2}+t)}{{{{\Gamma(2{{t}})\Gamma(-2{{t}})}}}%
\Gamma(z-\frac{\nu}{2}-t)\Gamma(z-\frac{\nu}{2}+t)}\\
&  \cdot{{\prod\limits_{k=1}^{4}\Gamma(-s_{k}-{\frac{{{\nu}}}{2}}%
+t)\Gamma(-s_{k}}}-{{{\frac{{{\nu}}}{2}}-t)dt}}\\
&  =\pi\int_{-\infty}^{\infty}\frac{\Gamma(1+z+\frac{\nu}{2}-ix)\Gamma
(1+z+\frac{\nu}{2}+ix)}{{{{\Gamma(2ix)\Gamma(-2ix)}}}\Gamma(z-\frac{\nu}%
{2}-ix)\Gamma(z-\frac{\nu}{2}+ix)}\\
&  \cdot{{\prod\limits_{k=1}^{4}\Gamma(-s_{k}-{\frac{{{\nu}}}{2}}%
+ix)\Gamma(-s_{k}}}-{{{\frac{{{\nu}}}{2}}-ix)dx}}.
\end{align*}
Using the integral representation given in \cite{rahman1986} as%
\begin{align}
&  \frac{1}{2\pi}\int_{-\infty}^{\infty}\frac{\Gamma(\lambda+ix)\Gamma
(\lambda-ix)\Gamma(\mu+ix)\Gamma(\mu-ix)}{{{{\Gamma(2i{x})\Gamma(-2{ix})}}}%
}\nonumber\\
&  \cdot\frac{\Gamma(\gamma+ix)\Gamma(\gamma-ix)\Gamma(\rho+ix)\Gamma
(\rho-ix)\Gamma(\sigma+ix)\Gamma(\sigma-ix)}{\Gamma(\tau+ix)\Gamma(\tau
-ix)}dx\nonumber\\
&  =\frac{2\Gamma(\lambda+\mu)\Gamma(\lambda+\gamma)\Gamma(\lambda+\rho
)\Gamma(\lambda+\sigma)}{\Gamma(\lambda+\tau)\Gamma(\mu+\tau)\Gamma
(\gamma+\tau)\Gamma(\lambda+\mu+\nu+\rho)}\nonumber\\
&  \cdot\frac{\Gamma(\mu+\gamma)\Gamma(\mu+\rho)\Gamma(\mu+\sigma
)\Gamma(\gamma+\rho)\Gamma(\gamma+\sigma)\Gamma(\lambda+\mu+\gamma+\tau
)}{\Gamma(\lambda+\mu+\nu+\sigma)}\nonumber\\
&  \cdot_{7}F_{6}%
\genfrac{[}{.}{0pt}{}{\lambda+\mu+\gamma+\tau-1,\frac{\lambda+\mu+\gamma
+\tau+1}{2},}{\frac{\lambda+\mu+\gamma+\tau-1}{2},\lambda+\tau,}%
\nonumber\\
&
\genfrac{.}{.}{0pt}{}{\lambda+\mu,\lambda+\gamma,\mu+\gamma,\tau-\sigma
,\tau-\rho}{\gamma+\tau,\mu+\tau,\lambda+\mu+\gamma+\sigma,\lambda+\mu
+\gamma+\rho}%
;1%
\genfrac{.}{]}{0pt}{}{{}}{{}}%
, \label{rah}%
\end{align}
where $_{7}F_{6}$ is the generalized hypergeometric series, and letting
$\lambda=1+z+\frac{\nu}{2},\mu=-s_{1}-\frac{\nu}{2},\gamma=-s_{2}-\frac{\nu
}{2},\rho=-s_{3}-\frac{\nu}{2},$and $\sigma=-s_{4}-\frac{\nu}{2},\tau
=z-\frac{\nu}{2}$ in Eq. (\ref{rah}), we simplify the solution $y_{\nu}(z)$
as
\begin{align*}
y_{\nu}(z)  &  =\frac{4\pi^{2}\Gamma(1+z-s_{1})\Gamma(1+z-s_{2})\Gamma
(1+z-s_{3})\Gamma(1+z-s_{4})\Gamma(-s_{1}-s_{2}-\nu)}{\Gamma(1+2z)\Gamma
(z-s_{1}-\nu)\Gamma(z-s_{2}-\nu)\Gamma(1+z-s_{1}-s_{2}-s_{3}-\nu)}\\
&  \cdot\frac{\Gamma(-s_{1}-s_{3}-\nu)\Gamma(-s_{1}-s_{4}-\nu)\Gamma
(-s_{2}-s_{3}-\nu)\Gamma(-s_{2}-s_{4}-\nu)}{\Gamma(1+z-s_{1}-s_{2}-s_{4}-\nu
)}\\
&  \cdot\Gamma(1+2z-s_{1}-s_{2}-\nu)\\
\cdot &  _{7}F_{6}%
\genfrac{[}{.}{0pt}{}{2z-s_{1}-s_{2}-\nu,\frac{2z-s_{1}-s_{2}-\nu+2}%
{2},}{\frac{2z-s_{1}-s_{2}-\nu}{2},1+2z,}%
\\
&
\genfrac{.}{.}{0pt}{}{1+z-s_{1},1+z-s_{2},-s_{1}-s_{2}-\nu,z+s_{4}%
,z+s_{3}}{z-s_{2}-\nu,z-s_{1}-\nu,1+z-s_{1}-s_{2}-s_{4}-\nu,1+z-s_{1}%
-s_{2}-s_{3}-\nu}%
;1%
\genfrac{.}{]}{0pt}{}{{}}{{}}%
,
\end{align*}
Ignoring a constant factor, the solution can be further written as%
\begin{align}
y_{\nu}(z)  &  =\frac{\Gamma(1+2z-s_{1}-s_{2}-\nu)\Gamma(1+z-s_{1}%
)\Gamma(1+z-s_{2})}{\Gamma(1+2z)\Gamma(z-s_{1}-\nu)\Gamma(z-s_{2}-\nu
)}\nonumber\\
&  \cdot\frac{\Gamma(1+z-s_{3})\Gamma(1+z-s_{4})}{\Gamma(1+z-s_{1}-s_{2}%
-s_{3}-\nu)\Gamma(1+z-s_{1}-s_{2}-s_{4}-\nu)}\nonumber\\
&  \cdot_{7}F_{6}%
\genfrac{[}{.}{0pt}{}{2z-s_{1}-s_{2}-\nu,\frac{2z-s_{1}-s_{2}-\nu+2}{2}%
,-s_{1}-s_{2}-\nu,}{\frac{2z-s_{1}-s_{2}-\nu}{2},1+2z,}%
\nonumber\\
&
\genfrac{.}{.}{0pt}{}{1+z-s_{1},1+z-s_{2},z+s_{4},z+s_{3}}{z-s_{2}-\nu
,z-s_{1}-\nu,1+z-s_{1}-s_{2}-s_{4}-\nu,1+z-s_{1}-s_{2}-s_{3}-\nu}%
;1%
\genfrac{.}{]}{0pt}{}{{}}{{}}%
\nonumber\\
&  =\frac{\prod\limits_{k=1}^{2}(z-s_{k}-\nu)_{\nu+1}\prod\limits_{k=3}%
^{4}(1+z-s_{1}-s_{2}-s_{k}-\nu)_{s_{1}+s_{2}+\nu}}{(1+2z-s_{1}-s_{2}%
-\nu)_{s_{1}+s_{2}+\nu}}\nonumber\\
&  \cdot_{7}F_{6}%
\genfrac{[}{.}{0pt}{}{2z-s_{1}-s_{2}-\nu,\frac{2z-s_{1}-s_{2}-\nu+2}{2}%
,-s_{1}-s_{2}-\nu,}{\frac{2z-s_{1}-s_{2}-\nu}{2},1+2z,}%
\nonumber\\
&
\genfrac{.}{.}{0pt}{}{1+z-s_{1},1+z-s_{2},z+s_{4},z+s_{3}}{z-s_{2}-\nu
,z-s_{1}-\nu,1+z-s_{1}-s_{2}-s_{4}-\nu,1+z-s_{1}-s_{2}-s_{3}-\nu}%
;1%
\genfrac{.}{]}{0pt}{}{{}}{{}}%
. \label{rah2}%
\end{align}

In particular, if we choose $s_{1}+s_{2}+\nu=n$ and $\nu+1=n$ in \ Eq.
(\ref{rah2}), then $s_{1}+s_{2}=1$. For this case, we can obtain a polynomial
solution for Eq. (\ref{th52}) as%
\begin{align}
y_{n}(z)  &  =\frac{\prod\limits_{k=1}^{4}(1+z-s_{k}-n)_{n}}{(2z+1-n)_{n}%
}\cdot_{7}F_{6}%
\genfrac{[}{.}{0pt}{}{2z-n,z-\frac{n}{2}+1,-n,}{z-\frac{n}{2},1+2z,}%
\nonumber\\
&
\genfrac{.}{.}{0pt}{}{z+s_{2},z+s_{1},z+s_{4},z+s_{3}}{1+z-s_{2}%
-n,1+z-s_{1}-n,1+z-s_{4}-n,1+z-s_{3}-n}%
;1%
\genfrac{.}{]}{0pt}{}{{}}{{}}%
,
\end{align}
which is the same as the well-known formula obtained in \cite{nikiforov1991},
(seen in Eq. (3.11.6) on page 134 in \cite{nikiforov1991}). This indicates
that Theorem 5.2 gives a more general solution form which includes the
well-known polynomial solution as its particular solution.

\section{{\protect\Large Conclusion }}

We have obtained the adjoint difference equation for the
Nikiforov-Uvarov-Suslov difference equation of hypergeometric type on
non-uniform lattices given as $x(s)=c_{1}q^{s}+c_{2}q^{-s}+c_{3}$ or
$x(s)=\tilde{c}_{1}s^{2}+\tilde{c}_{2}s+\tilde{c}_{3}$ and proved it to be a
difference equation of hypergeometric type on non-uniform lattices as well.
The particular solutions of the adjoint difference equation have then been
obtained. By applying these particular solutions for the adjoint equation, we
can obtain the particular solutions of the original difference equation of
hypergeometric type on non-uniform lattices. Finally, we have obtained new
fundamental theorems for the Nikiforov-Uvarov-Suslov difference equation of
hypergeometric type and illustrated their applications by an example.

\vspace{0.2in}

\hspace{-0.2in}\textbf{Acknowledgements. }The first author was supported by
the Fundamental Research Funds for the Central Universities of China, grant
number 20720150006, and Natural Science Foundation of Fujian province of
China, grant number 2016J01032.

\end{document}